\documentclass[11pt]{IEEEtran}
\pdfoutput=1
\usepackage[utf8]{inputenc}
\usepackage[english]{babel}
\fontfamily{phv}\selectfont
\usepackage{stfloats}
\usepackage{cite}
\usepackage{amsthm,amsmath,amssymb,amsfonts,mathtools}
\usepackage{tabularx}
\usepackage[caption=false,font=footnotesize]{subfig}
\usepackage{graphicx}
\usepackage{textcomp}
\usepackage[hyphens]{url}
\usepackage[hidelinks,breaklinks]{hyperref}
\usepackage[normalem]{ulem}

\begin{document}
\title{Benchmarking the performance of controllers for power grid transient stability}
\author{Randall~Martyr,
	Benjamin~Sch\"{a}fer,
	Christian~Beck,
	and~Vito~Latora
	\thanks{Financial support received from the UK Engineering and Physical Sciences Research Council (EPSRC Reference: EP/N013492/1) and the German Federal Ministry of Education and Research (BMBF grant no. 03SF0472A-F).}
	\thanks{Randall Martyr, Christian Beck and Vito Latora are with the School of Mathematical Sciences, Queen Mary University of London, Mile End Road, London E1 4NS, United Kingdom (e-mail: r.martyr@qmul.ac.uk; c.beck@qmul.ac.uk; v.latora@qmul.ac.uk).}
	\thanks{Vito Latora is also with the Dipartimento di Fisica ed Astronomia, Universit\`a di Catania and INFN, I-95123 Catania, Italy.}
	\thanks{Benjamin Sch\"{a}fer is with the Chair for Network Dynamics, Technical University of Dresden and Max Planck Institute for Dynamics and Self-Organization (MPIDS), 37077 G\"{o}ttingen, Germany (e-mail: benjamin.schaefer@ds.mpg.de).}
}

\maketitle

\begin{abstract}
As the energy transition transforms power grids across the globe, it poses several challenges regarding grid design and control. In particular, high levels of intermittent renewable generation complicate the task of continuously balancing power supply and demand, requiring sufficient control actions. Although there exist several proposals to control the grid, most of them have not demonstrated to be cost efficient in terms of optimal control theory. Here, we mathematically formulate an optimal {\em centralized} (therefore non-local) control problem for stable operation of power grids and determine the minimal amount of active power necessary to guarantee a stable service within the operational constraints, minimizing a suitable cost function at the same time. This optimal control can be used to benchmark control proposals and we demonstrate this benchmarking process by investigating the performance of three distributed controllers, two of which are fully decentralized, that have been recently studied in the physics and power systems engineering literature. Our results show that cost efficient controllers distribute the controlled response amongst all nodes in the power grid. Additionally, superior performance can be achieved by incorporating sufficient information about the disturbance causing the instability. Overall, our results can help design and benchmark secure and cost-efficient controllers.
\end{abstract}


\section{Introduction}

\begin{table}
	\renewcommand{\arraystretch}{1.3}
	\caption{Nomenclature used in this paper. Vectors and matrices are denoted in boldface.}
	\label{Table:Parameters}
	\centering
	\begin{tabularx}{\linewidth}{|cXc|}
		\hline
		{\bf Notation} & {\bf Description} & {\bf Units} \\[5pt]
		& \centering {\bf Synchronous machine parameters} & \\
		$\mathcal{N}$ & the set of nodes $\{1,\ldots,N\}$ where $N \ge 2$ is the number of nodes in the network & -- \\
		$\Omega$ & synchronous angular velocity used as reference & $rad \cdot s^{-1}$ \\
		$\boldsymbol{B}$ & $N \times N$-dimensional matrix of line susceptances & pu \\
		$M_{i}$ & inertia coefficient & $s^{2}$ \\
		$D_{i}$ & damping coefficient & pu \\
		$E_{f,i}$ & exciter voltage & pu \\
		$X_{d,i}$ & direct synchronous reactance & pu \\
		$X'_{d,i}$ & direct synchronous transient reactance & pu \\
		$T'_{do,i}$ & direct axis transient time constant & $s$ \\
		$P_{e,i}$ & electromagnetic air-gap power & pu \\
		$P_{in,i}$ & net power injection, the difference between mechanical power and aggregate load & pu \\[5pt]
		& \centering {\bf Synchronous machine state quantities} & \\
		$\theta_{i}$ & rotor angle relative to the grid reference & $rad$ \\
		$\omega_{i}$ & angular velocity relative to the grid reference & $rad \cdot s^{-1}$ \\
		$V_{i}$ & normalized machine voltage & pu \\
		$\xi_{i}$ & disturbance to net power injection & pu \\
		$\sigma(\boldsymbol{\omega})$ & standard deviation of network angular velocities & $rad \cdot s^{-1}$ \\
		$\langle \omega \rangle$ & mean value of network angular velocities & $rad \cdot s^{-1}$ \\[5pt]
		& \centering {\bf Optimization parameters and variables} & \\
		$T$ & control time horizon & $s$ \\
		$\boldsymbol{x}$ & $3N$-dimensional state vector & -- \\
		$\boldsymbol{u}$ & $N$-dimensional vector of controlled power injections & pu \\
		$\mathcal{U}$ & set of control variables & -- \\
		$J$, $C_{\eta}$, $\varepsilon_{\eta}$ & cost functional, constraint functional, constraint tolerance & --  \\\hline
	\end{tabularx}
\end{table}

The electrical power grid is undergoing drastic changes due to the energy transition \cite{Turner1999,Boyle2004,Ueckerdt2015} and suitable control approaches are necessary to ensure a reliable and stable operation \cite{Kundur1994}. The generation side of the grid is changing as additional renewable generators are installed to mitigate climate change, introducing fluctuations on a time scale of days \cite{Heide2010} to sub-seconds \cite{Milan2013}. In addition, the demand side is changing due to the ongoing electrification of heating and transport \cite{Dennis2016} and the introduction of {\it demand control} \cite{Palensky2011}. Regardless of these changing conditions, the grid needs to stay within strict operational boundaries to guarantee a stable electricity supply and to prevent damage to sensitive electronic devices \cite{Kundur1994}.

A fundamental aspect of power system stability is the ability of interconnected synchronous machines of a power system to remain synchronized. {\it Transient stability} describes the power system's ability to maintain synchronism in the face of severe transient disturbances \cite{Kundur1994}, and is of great importance in preventing cascading failures \cite{Motter2002,Yang2017,Schaefer2017}. Control mechanisms that balance active power and regulate frequency in the grid are key to maintaining these stability conditions. {\it Primary controls} \cite{Zhao2014} respond within a few seconds of an event to stabilize the frequency within its permissible operating limits, after which {\it secondary} \cite{Xi2017,Tchuisseu2018} and {\it tertiary controls} restore the frequency to its nominal value \cite{Kirschen2004}.

In this paper we describe control algorithms for networked systems (such as the power grid) as being {\it centralized} if a central controller performs computations and issues control actions for the entire network, {\it distributed} if there are multiple autonomous controllers that perform computations and can communicate with each other, and {\it decentralized} if there are multiple autonomous controllers that perform computations but do not communicate with each other. Our definition intentionally permits distributed controllers that do not communicate with each other, thus making decentralized controllers a special case, albeit degenerate. Distributed approaches are often supported via advanced power electronics \cite{Carrasco2006} and economic considerations \cite{Li2016} to further improve the grid's stability. For large-scale networks, centralized control schemes can be computationally complex and impractical, thereby making distributed control schemes with low computation and communication complexity more desirable \cite{Molzahn2017}. Decentralized controllers are popular choices since they rely only on local measurements, but they can have poor system-wide performance in practice \cite{Molzahn2017,Venkat2008}. For a discussion on the strengths and limitations of centralized, decentralized and distributed controllers for power systems see \cite{Molzahn2017}.

In this paper we seek to answer the following question: What are the characteristics of a controller that efficiently synchronizes the power grid in the presence of known disturbances caused by changes in demand and generation? We answer this question by investigating the solution to an optimal control problem (see \cite{Hestenes1966,Clarke2013}) for synchronization of a power grid described by a network of control areas (nodes) $\mathcal{N}$. Note that the optimal control has complete information regarding the temporal evolution of the disturbance at all nodes in the network. Therefore, it constitutes the ideal controller in terms of performance and any realistic controller, centralized or distributed, can be compared in its performance to the optimal one. In this paper, we use the optimal control to exemplarily benchmark the following three distributed control schemes, two of which are fully decentralized.

Sch\"{a}fer et al \cite{Schaefer2015,Schaefer2016} recently investigated a decentralized {\it linear local frequency} (LLF) controller, linked to a patent \cite{Walter2016}, that can improve the grid's transient stability by regulating electricity demand and supply through economic incentives. The control action at area $i \in \mathcal{N}$ is directly proportional to $\omega_{i}$, the local angular velocity deviation relative to the grid reference,
\begin{equation}\label{eq:Control-Proportional-Frequency}
u_{i}(t) \coloneqq -\nu_{i}\omega_{i}(t) \quad i \in \mathcal{N},
\end{equation}
with $\nu_{i} > 0$. The constant $\nu_{i}$ in \eqref{eq:Control-Proportional-Frequency} measures the willingness at node $i$ to change the active power level and effectively increases the damping parameter from $D_{i}$ to $D_{i} + \nu_{i}$ in the grid dynamics \eqref{eq:3rd-order-controlled-swing-equation} below.

In \cite{Zhao2015,Weitenberg2017} the following {\it integral local frequency} (ILF) control is studied,
\begin{equation}\label{eq:Control-Integral-Frequency}
u_{i}(t) \coloneqq -\frac{1}{\kappa_{i}}\int_{0}^{t}\omega_{i}(\tau){d}\tau, \quad i \in \mathcal{N},
\end{equation}
where $\kappa_{i} > 0$. The integral control \eqref{eq:Control-Integral-Frequency} can improve the power grid's synchronization and stability, and can be economically efficient in a particular sense \cite{Zhao2015,Weitenberg2017}.

Finally, we consider the following gather-and-broadcast (GAB) distributed controller which is a special case of the one defined in \cite{Doerfler2017},
\begin{equation}\label{eq:GAB-Control-Local}
u_{i}(t) \coloneqq -\frac{1}{\mu_{i}}\int_{0}^{t}\sum_{j=1}^{N}A_{ij}\omega_{j}(\tau){d}\tau, \quad i \in \mathcal{N},
\end{equation}
where $\mu_{i} > 0$ and $\mathbf{A} = (A_{ij})_{(i,j) \in \mathcal{N} \times \mathcal{N}}$ is an unweighted adjacency matrix, $A_{ij} \in \{0,1\}$ and $A_{ij} = A_{ji}$, that defines a communication network between the control areas. If $A_{ij} = 1$ when $i = j$ and $A_{ij} = 0$ otherwise, then the GAB controller \eqref{eq:GAB-Control-Local} reduces to the decentralized integral controller \eqref{eq:Control-Integral-Frequency}. In this paper we consider the special case of a fully connected communication network, $A_{ij} = 1$ for all $(i,j) \in \mathcal{N} \times \mathcal{N}$, which leads to,
\begin{align}\label{eq:GAB-Control-Global} 
u_{i}(t) & \coloneqq -\frac{1}{\mu_{i}}\int_{0}^{t}\sum_{j=1}^{N}\omega_{j}(\tau){d}\tau \nonumber \\
& = -\frac{N}{\mu_{i}}\int_{0}^{t}\frac{\sum_{j=1}^{N}\omega_{j}(\tau)}{N}{d}\tau, \quad i \in \mathcal{N},
\end{align}
thereby making the GAB controller proportional to the time integral of the mean angular velocity.

In the following section we present the optimal control problem for power grid synchronization. The power grid dynamics are given by a system of ordinary differential equations for a state vector $\boldsymbol{x}$ of phase angles, angular velocity deviations (related to the grid frequency) and voltage amplitudes. Let $\mathcal{U}$ be a suitable set of time-dependent control variables $\boldsymbol{u}$. For a given $\boldsymbol{u} \in \mathcal{U}$, we quantify its cost through a cost function $J(\boldsymbol{u})$, and evaluate its performance with respect to various operational constraints $C_{\eta}(\boldsymbol{u})$ and their tolerances $\varepsilon_{\eta}$. The optimal control problem for power grid synchronization is expressed mathematically as follows.
\begin{quotation}
	\noindent{\bf Problem:}
	\begin{gather}
	\text{minimize}\; J(\boldsymbol{u})\; \text{subject to:} \nonumber \\
	\begin{split}
	i) & \quad
	\dot{\boldsymbol{x}}(t)  = \boldsymbol{f}(t,\boldsymbol{x}(t),\boldsymbol{u}(t)), \quad
	\boldsymbol{x}(0) = \boldsymbol{x}_{0}; \\
	ii) & \quad \boldsymbol{u} \in \mathcal{U}; \\
	iii) & \quad C_{\eta}(\boldsymbol{u}) \le \varepsilon_{\eta}  \enskip \text{for} \enskip \eta = 1,\ldots,N+2,
	\end{split}\label{Problem:Optimal-Control}
	\end{gather}
	where $\boldsymbol{f}$ governs the intrinsic dynamics of the state of the grid (see \eqref{eq:3rd-order-controlled-swing-equation} below), and $N \ge 2$ is the number of nodes in its representation as a network.
\end{quotation}
Problem~\eqref{Problem:Optimal-Control} is solved numerically using a control parametrization method \cite{Teo1991} that is outlined in the Appendix. In Section \ref{Section:Simulations} we illustrate the efficiency of the optimal control compared to the three proposed controls, \eqref{eq:Control-Proportional-Frequency}, \eqref{eq:Control-Integral-Frequency} and \eqref{eq:GAB-Control-Global}, for a four-node network motif. Finally, in Section \ref{Section:Outlook} we close with a conclusion and outlook.

Our results show that the optimal control achieves superior performance with respect to cost whilst achieving comparable and, in some respects, better performance with respect to the operational constraints. However, this superiority is a consequence of the optimal control utilizing its knowledge of the disturbance to form a pre-emptive response. Realistic controllers will not have this information for random disturbances and will therefore require larger investments than the optimal control. Nevertheless, since the distributed controllers we investigate do not explicitly incorporate any information about the disturbance, we postulate that realistic controllers can achieve superior performance if they incorporate some of this information. Regularly occurring disturbances, for instance those caused by economic effects \cite{schaefer2018} or steep gradients due to the sun rising (similar to the recent solar eclipse) \cite{harrison2016solar}, provide important examples in which information about the disturbance may be obtained practically.

\section{An optimal control problem for power grid transient stability}\label{Section:Optimal-Control}
This section details the optimal control problem \eqref{Problem:Optimal-Control} that we use to benchmark the distributed (including decentralized) controllers' performances. However, before focusing on optimal control we need to discuss the model that we use for the intrinsic dynamics of the power grid.

\subsection{Dynamics for transient stability analysis}\label{Section:Dynamics}
The rotor mechanical velocities of the interconnected synchronous machines in a power grid must be synchronized to the same frequency, else there can be deviations in the rotor angles that lead to instabilities \cite[p.~19]{Kundur1994}. A severe transient disturbance can cause large deviations in the rotor angles, which may lead to a progressive drop in the nodal voltages \cite[p.~27]{Kundur1994} and further affect the angular velocities and rotor angle values. A realistic model of the power grid should therefore take the influence of the rotor angles' deviations on the voltage amplitudes into account. This allows us to analyze slower phenomena such as large deviations in voltage or frequency, as typically done in mid-term stability studies \cite[p.~34]{Kundur1994}. Therefore, in this paper we use a third order model \cite[p.~456]{Machowski2008}, which describes the power grid as a network of $N \ge 2$ control areas, each represented by a synchronous generator or motor and governed by a set of differential equations for the rotor angle $\theta_i$, angular velocity deviation $\omega_i$, and voltage $V_i$ at each node,
\begin{equation}
\begin{split}
(i) &  \quad \dot{\theta}_{i} = \omega_{i} \\
(ii) &  \quad
M_{i}\dot{\omega}_{i} = P_{in,i} - P_{e,i} + u_{i} - D_{i}\omega_{i} \\
(iii) &  \quad  T'_{do,i}\dot{V}_{i} = E_{f,i}-V_{i} + I_{d,i}\bigl(X_{d,i} - X'_{d,i}\bigr),\\
\text{for}& \; i=1,2,\ldots,N,
\end{split}\label{eq:3rd-order-controlled-swing-equation}
\end{equation}
where $P_{in,i}$ is the net power injection, $u_{i}$ is the controlled active power, $P_{e,i}$ is the electrical power,
\[
P_{e,i} = \sum_{j = 1}^{N}B_{i,j}\sin(\theta_{i}-\theta_{j})V_{i}V_{j},
\]
$I_{d,i}$ is the armature current,
\[
I_{d,i} = \sum_{j = 1}^{N}B_{i,j}\cos(\theta_{i}-\theta_{j})V_{j},
\]
and $M_{i}$, $D_{i}$, $T'_{do,i}$, $X_{d,i}$, $X'_{d,i}$ and $E_{f,i}$ are parameters described in Table~\ref{Table:Parameters}. This model assumes a lossless network and a constant exciter voltage (emf) $E_{f,i}$. It also neglects transient saliency power and ignores damping effects produced by eddy currents. Note that $\omega_i$ represents the deviation of the rotor angle velocity from a synchronized state $2\pi F$, where $F$ is the reference frequency in Hertz. However, for brevity we will often say ``angular velocity'' instead of ``angular velocity deviation''.

A positive value for $P_{in,i}$ indicates net generation at node $i$ and in this case we refer to this node as a generator. A negative value of $P_{in,i}$ indicates net consumption at node $i$ and in this case we refer to this node as a consumer or motor. We refer to positive values for the control variable $u_{i}$ as {\it incremental actions} \cite{Szabo2017} since they correspond to an increase in generation or an equivalent decrease in demand. Similarly, we refer to negative values for $u_{i}$ as {\it decremental actions} \cite{Szabo2017} since they correspond to a decrease in generation or an equivalent increase in demand.

\subsection{Operational constraints of the power grid}\label{Section:Constraints}
Let $\boldsymbol{x} = (x_{1},\ldots,x_{3N})$ denote the $3N$-dimensional controlled state variable obtained from \eqref{eq:3rd-order-controlled-swing-equation} with components given by
\begin{equation}\label{eq:State-Variable-Components}
x_{i} = \theta_{i},\; x_{N+i} = \omega_{i},\; x_{2N+i} = V_{i} \enskip \text{for} \enskip i \in \mathcal{N}.
\end{equation}
The dynamics of $\boldsymbol{x}$ in \eqref{eq:3rd-order-controlled-swing-equation}
can be written compactly as
\begin{equation}\label{eq:Equations-Of-Motion}
\dot{\boldsymbol{x}}(t) = \boldsymbol{f}(t,\boldsymbol{x}(t),\boldsymbol{u}(t)),
\end{equation}
where expressions for the components of the intrinsic dynamics $\boldsymbol{f} = (f_{1},\ldots,f_{3N})$ are obtained from \eqref{eq:3rd-order-controlled-swing-equation} using the assignment given in \eqref{eq:State-Variable-Components}. Each component of the control variable $\boldsymbol{u} = (u_{1},\ldots,u_{n})$ corresponds to the amount of additional active power injected or withdrawn at an individual node in the network. We assume that controls are bounded: for each $i \in \mathcal{N}$ we have $u_{i}(t) \in U_{i}$ where:
\begin{equation}\label{eq:Control-Set}
U_{i} = [u^{min}_{i},u^{max}_{i}], \quad -\infty < u^{min}_{i} < u^{max}_{i} < \infty.
\end{equation}
Let $\mathcal{U}$ denote the set of all such control functions.

\paragraph*{\bf Synchronization}
In our model, synchronization of the rotor angle velocities for the control areas means $\omega_{i} = \omega_{j}$ for all $i,j \in \mathcal{N}$. Letting $\boldsymbol{\omega} = (\omega_{1},\ldots,\omega_{N})$ denote the vector of angular velocities and $\langle \omega \rangle = \frac{1}{N}\sum_{j = 1}^{N}\omega_{j}$ its arithmetic mean, we measure the lack of synchronization using the standard deviation of $\boldsymbol{\omega}$,
\begin{equation}\label{eq:Effective-Frequency-Dispersion}
\sigma(\boldsymbol{\omega}) = \left(\frac{1}{N}\sum_{i=1}^{N}\left(\omega_{i} - \langle \omega \rangle\right)^{2}\right)^{\frac{1}{2}}.
\end{equation}
Let $0 < T < \infty$ denote the length of the control horizon $[0,T]$ in seconds. Define the synchronization constraint loss function by
\begin{equation}\label{eq:Stability-Constraint}
\psi_{1}(\boldsymbol{x}) = -\sigma(\boldsymbol{\omega}),
\end{equation}
and the total synchronization loss on $[0,T]$ by
\begin{align}
C_{1}(\boldsymbol{u}) = {} & \int_{0}^{T}\left(\min(0,\psi_{1}(\boldsymbol{x}(t)))\right)^{2}{d}t \nonumber \\
& + \lambda_{1}\min\bigl(0,\psi_{1}(\boldsymbol{x}(T)\bigr)^{2} \nonumber \\
= {} & \int_{0}^{T}\sigma(\boldsymbol{\omega}(t))^{2}{d}t  + \lambda_{1}\sigma(\boldsymbol{\omega}(T))^{2},\label{eq:Synchronization-Total-Loss}
\end{align}
where $\lambda_{1} \ge 0$ is a weight parameter which emphasizes the relative importance of the constraint at the final time $T$. Recalling the definition of $\sigma(\boldsymbol{\omega})$ in \eqref{eq:Effective-Frequency-Dispersion}, the quadratic weighting given to it naturally defines the variance of $\boldsymbol{\omega}$. Other weighting schemes are also possible.

\paragraph*{\bf Mean angular velocity operational limits}
The variable $\omega_{i}$ quantifies the deviation of the angular velocity at node $i$ from the synchronous reference $\Omega$ (rad/s), where $\Omega$ is related to the nominal frequency $F$ (Hz) of the power grid by $\Omega = 2 \pi \cdot F$. In the United Kingdom and many other countries the nominal frequency is $F = 50$ Hz. For reasons related to the quality of electricity supply, the frequency must respect certain operational limits. In the United Kingdom, for example, the statutory limits are $\pm$0.5 Hz of the nominal value 50 Hz, and the operational limits are set to the stricter range of $\pm$0.2 Hz \cite{nationalgrid1}. In our model, this implies the values of the mean angular velocity $\langle \omega \rangle$ should be constrained,
\begin{equation}\label{eq:Mean-Angular-Velocity-State-Constraint}
\omega_{min} \le \langle \omega \rangle \le \omega_{max}.
\end{equation}
Define the mean angular velocity constraint loss function by
\begin{equation}\label{eq:Angular-Velocity-Constraints}
\psi_{2}(\boldsymbol{x}) = (\omega_{max} - \langle \omega \rangle)(\langle \omega \rangle - \omega_{min}),
\end{equation}
and the total loss on $[0,T]$ for violating this constraint by
\[
\begin{split}
C_{2}(\boldsymbol{u}) = {} & \int_{0}^{T}\left(\min(0,\psi_{2}(\boldsymbol{x}(t)))\right)^{2}{d}t \\
& + \lambda_{2}\min\bigl(0,\psi_{2}(\boldsymbol{x}(T))\bigr)^{2},
\end{split}
\]
where $\lambda_{2} \ge 0$ is a weight parameter. Note that only when $\psi_{2}(\boldsymbol{x})$ is
negative in eq.~(\ref{eq:Angular-Velocity-Constraints}) we get a contribution.
\paragraph*{\bf Voltage operational limits}
Since the voltages in our model are also time dependent, it is important to also take into account appropriate operational constraints on these variables. For example, regulations in the United Kingdom require that the steady state voltages should be kept within $\pm 6\%$ of the nominal voltage for systems between 1 and 132 (kV), or $\pm 10\%$ of the nominal voltage for systems above 132 (kV) \cite{nationalarchives}. In our model we can take this into account with the following constraint,
\begin{equation}\label{eq:Voltage-Constraint}
V_{i}^{min} \le V_{i} \le V_{i}^{max},\enskip i \in \mathcal{N},
\end{equation}
where $V_{i}^{min} < V_{i}^{max}$. We define a loss function for the voltage constraint at each node $i \in \mathcal{N}$ by
\begin{equation}\label{eq:Voltage-Constraints}
\psi_{2+i}(\boldsymbol{x}) = (V_{i}^{max} - V_{i})(V_{i} - V_{i}^{min}),
\end{equation}
and the total loss on $[0,T]$ for violating this constraint by
\[
\begin{split}
C_{2+i}(\boldsymbol{u}) = {} & \int_{0}^{T}\left(\min(0,\psi_{2+i}(\boldsymbol{x}(t)))\right)^{2}{d}t \\
& + \lambda_{2+i}\min\bigl(0,\psi_{2+i}(\boldsymbol{x}(T))\bigr)^{2},
\end{split}
\]
where $\lambda_{2+i} \ge 0$ are weight parameters.

\subsection{Formulation of the optimal control problem}
For $\eta = 1,\ldots,N+2$ the total loss $C_{\eta}$ is non-negative, and is equal to zero if, equivalently, the $\eta$-th constraint is satisfied on $[0,T]$. We relax this by introducing {\it tolerance parameters} $\varepsilon_{\eta} \ge 0$, $\eta \in \{1,\ldots,N+2\}$, and say that a control $\boldsymbol{u} \in \mathcal{U}$ is {\it feasible} if it satisfies
\begin{equation}\label{eq:Canonical-Inequality-Constraints}
C_{\eta}(\boldsymbol{u}) \le \varepsilon_{\eta}  \enskip \text{for} \enskip \eta = 1,\ldots,N+2.
\end{equation}
Below we define a cost objective $J(\boldsymbol{u})$ which we use with the constraint losses \eqref{eq:Canonical-Inequality-Constraints} to formulate the optimal control problem \eqref{Problem:Optimal-Control}.

At an initial time $t = 0$, the power grid is synchronized and at a steady state, $\dot{\boldsymbol{x}} = 0$, in which various operational constraints are satisfied. Suppose the constant net injection $P_{in,i}$ corresponding to the steady state is perturbed according to an external disturbance $\xi_{i}$,
\[
P_{in,i} \to P_{in,i} + \xi_{i}(t), \quad t \in [0,T],
\]
that causes the grid to become unsynchronized. We would like the control function $\boldsymbol{u}$ to return the grid close to a synchronized state before $T$ seconds, and with a ``minimal cost'' that ensures the constraint conditions \eqref{eq:Canonical-Inequality-Constraints} are satisfied. Let $L(t,\boldsymbol{x},\boldsymbol{u})$ denote the value of a {\it cost rate function} $L$ that can generally depend on time and the current value of the state and control vectors. Letting $\boldsymbol{I}_{N}$ denote the $N \times N$ identity matrix and ${\rm tr}$ denote the transpose operator, we define the following quadratic cost,
\begin{equation}\label{eq:Euclidean-Control-Cost}
L(t,\boldsymbol{x},\boldsymbol{u}) \coloneqq \boldsymbol{u}^{{\rm tr}}\boldsymbol{I}_{N}\boldsymbol{u} = \sum_{i=1}^{N}(u_{i})^{2},
\end{equation}
which is typical of those in the frequency control literature \cite{Trip2016,Weitenberg2017}. The rate function \eqref{eq:Euclidean-Control-Cost} is used to define the following total cost for a control $\boldsymbol{u} \in \mathcal{U}$,
\begin{equation}
J(\boldsymbol{u}) = \int_{0}^{T}L\bigl(t,\boldsymbol{x}(t),\boldsymbol{u}(t)\bigr){d}t.\label{eq:Optimal-Control-Lagrange-Criterion}
\end{equation}

By the definition \eqref{eq:Euclidean-Control-Cost} of the cost rate, the cost objective \eqref{eq:Optimal-Control-Lagrange-Criterion} assigns higher costs to control functions $\boldsymbol{u}$ that exert large amounts of effort over time. Moreover, adjustments in demand and generation of the same magnitude are penalized equally due to the symmetry $\boldsymbol{u}^{{\rm tr}}\boldsymbol{I}_{N}\boldsymbol{u} = (-\boldsymbol{u})^{{\rm tr}}\boldsymbol{I}_{N}(-\boldsymbol{u})$. If demand and generation should be penalized differently then this can be achieved by adjusting \eqref{eq:Euclidean-Control-Cost}. Note that by using the identity matrix $\boldsymbol{I}_{N}$ we assume that the cost of control is independent of the node. If this is not the case, then we can replace $\boldsymbol{I}_{N}$ in \eqref{eq:Euclidean-Control-Cost} with another positive diagonal matrix. Finally, if we should also ensure that the system state does not deviate too far from its initial value $\boldsymbol{x}_{0}$, then we can penalize such deviations by adjusting the cost rate \eqref{eq:Euclidean-Control-Cost} or constraints $C_{\eta}$.

\section{Simulations for a four-node networked power system}\label{Section:Simulations}

\begin{figure*}[!th]
	\centering
	\subfloat[Four-node network motif]{\includegraphics[height=0.25\textheight, width=0.4\textwidth]{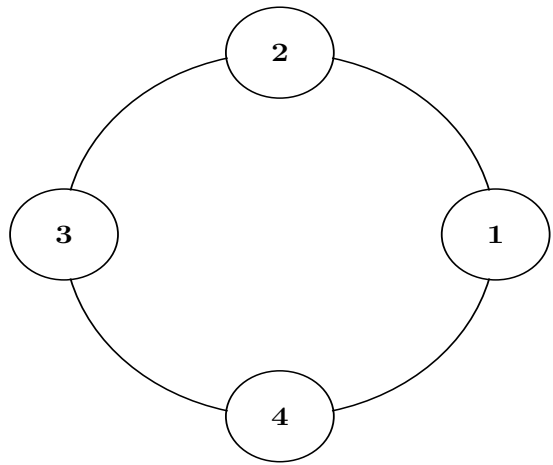}\label{fig:Network-4-Nodes}
	}
	\hfil
	\subfloat[Power disturbance $\xi_{1}$ to node 1]{\includegraphics[height=0.25\textheight, width=0.475\textwidth]{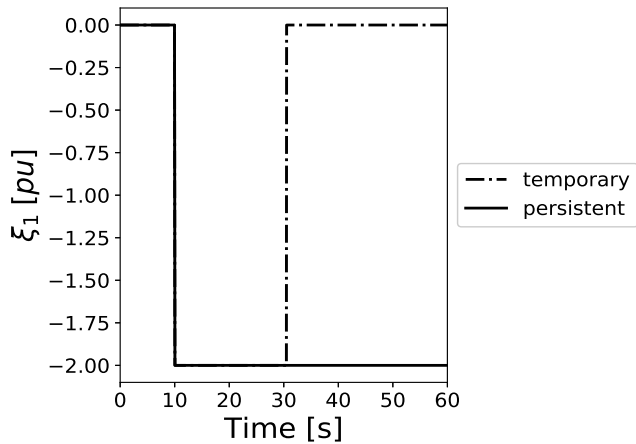}\label{fig:Mechanical-Power_Disturbances}
	}
	\caption{The first figure (a) illustrates the four-node motif network with ring topology. Parameters are given in the Appendix. The second figure (b) illustrates the types of power disturbance $\xi_{1}$ applied at node 1. We consider a short temporary change (dash-dotted line) of power and a persistent change (solid line). No disturbances are applied to the other nodes.}
\end{figure*}

For the numerical simulations we use the test system shown in Fig.~\ref{fig:Network-4-Nodes}. Note that such a network may be obtained as a reduction of a larger network, for example the IEEE 39-bus test system \cite{Ourari2006,Nabavi2013}. We consider two types of disturbance with each one altering the net power injection at node 1 as shown in Fig.~\ref{fig:Mechanical-Power_Disturbances}. The {\it temporary} disturbance reflects a sudden but short doubling of demand, or equivalent loss of generation, at node 1 from time $t = 10~\text{s}$ that lasts for only twenty seconds. The {\it persistent} disturbance reflects a sudden doubling of demand at node 1 from time $t = 10~\text{s}$ that lasts for the remaining control horizon. Results for the case with an analogous increase in generation, or equivalent loss of demand, are symmetric and thus omitted. In Appendix \ref{Section:Parameter-Tables} we list the parameter values for the model and control problem.

Upon representing the constraints by an appropriately defined vector of auxiliary state variables, we can apply the theoretical results in \cite{Clarke2013} or \cite{Cesari1983} to assert the existence of a solution to the optimal control problem \eqref{Problem:Optimal-Control}. Furthermore, Pontryagin's Maximum Principle \cite{Hestenes1966,Clarke2013} provides us with a set of mathematical conditions that a solution to the optimal control problem necessarily satisfies. Instead of pursuing this mathematical formalism, however, we empirically investigate characteristics of an optimal control by solving the optimal control problem numerically. The numerical solutions are obtained using the control parametrization method \cite{Teo1991}, which approximates the optimal control problem \eqref{Problem:Optimal-Control} by a constrained non-linear optimization problem over a bounded $(N \times n_{p})$-dimensional space, where $n_{p}$ is a positive integer, that parametrizes step control functions as follows,
\begin{equation}\label{eq:Step-Function-Control}
u_{i}(t) = \sum_{k = 1}^{n_{p}}u_{i}^{k}\mathbf{1}_{[t_{k-1},t_{k})}(t),\quad u_{i}^{k} \in U_{i},\; i \in \mathcal{N}.
\end{equation}
Further details of the algorithm are given in Appendix \ref{Appendix:Control-Parametrization-Method}, and the source code for the numerical experiments is available online \cite{Martyr2018}. For the simulations we use equidistant partitioning points $t_{k} = \frac{k}{n_{p}}T$, $0 \le k \le n_{p}$, with $n_{p} = 1500$, and the Sequential Least Squares Programming (SLSQP) routine in Python to solve the non-linear optimization problem.

We compare the performance of the optimal control (OC) and three controllers, LLF \eqref{eq:Control-Proportional-Frequency}, ILF \eqref{eq:Control-Integral-Frequency} and GAB \eqref{eq:GAB-Control-Global}, restricting values of the latter controls to the set $U = \prod_{i \in \mathcal{N}}U_{i}$ if necessary. We use the trapezoidal rule to approximate the integrals in \eqref{eq:Control-Integral-Frequency} and \eqref{eq:GAB-Control-Global} and update the control $u_{i}$ incrementally in an online manner. For simplicity we suppose that $\nu_{i} = \nu$, $\kappa_{i} = \kappa$ and $\mu_{i} = \mu$ for all $i \in \mathcal{N}$ in \eqref{eq:Control-Proportional-Frequency} and \eqref{eq:Control-Integral-Frequency} respectively.

The proposed distributed controllers are designed to keep the system frequency close to the nominal value (and, therefore, the angular velocity close to 0). In order to make the comparison fair we therefore choose the angular velocity constraints in \eqref{eq:Mean-Angular-Velocity-State-Constraint} to reflect a maximum allowed deviation of $0.1\%$ from the nominal value $50$ Hz, which is $\pm0.05$ Hz. Table~\ref{Table:Control-Coefficient-Values} below shows the values of $\nu$, $\kappa$ and $\mu$ we used in the simulations.

\begin{table}[!h]
	\renewcommand{\arraystretch}{1.3}
	\caption{Parameter values $\nu~(s^{-1})$, $\kappa~(s^{-2})$ and $\mu~(s^{-2})$ selected for the proposed controllers. Also included is the synchronization constraint total loss $C_{1}$ for the temporary (T) and persistent (P) disturbances and the respective performance of the optimal control (OC).}
	\label{Table:Control-Coefficient-Values}
	\centering
	\begin{tabular}{ccc}
		\hline
		{\bf Control} & {\bf Parameter Value} & {\bf $C_{1}$ (T, P)} \\ \hline
		LLF \eqref{eq:Control-Proportional-Frequency} & $\nu = 1$ & $3.6 \cdot 10^{-3}$,\, $1.8 \cdot 10^{-3}$ \\
		ILF \eqref{eq:Control-Integral-Frequency} & $\kappa = 15$ & $6.3 \cdot 10^{-3}$,\, $3 \cdot 10^{-3}$ \\
		GAB \eqref{eq:GAB-Control-Global} & $\mu = 60$ & $6.3 \cdot 10^{-3}$,\, $3 \cdot 10^{-3}$ \\
		OC \eqref{Problem:Optimal-Control} & -- & $10^{-4}$,\, $10^{-4}$
		\\\hline
	\end{tabular}
\end{table}

The value for $\nu$ was chosen to be comparable to the damping constants given in the Appendix. The value for $\mu$ was selected according to the simulations in \cite[p.~303]{Doerfler2017}, whilst the value for $\kappa$ was selected to satisfy $\frac{\mu}{\kappa} = N = 4$, based on the relation in \eqref{eq:GAB-Control-Global} above. Notice that the synchronization total loss $C_{1}$ \eqref{eq:Synchronization-Total-Loss} for the distributed controls is larger for the temporary disturbance than for the persistent one. This is because the temporary disturbance causes two sudden changes to the net power injection over the control horizon whereas the persistent disturbance only causes one sudden change.

\subsection{Simulated dynamics of the controlled power system}

\begin{figure*}[!th]
	\centering
	\subfloat[LLF for temporary disturbance]{
		\includegraphics[height=0.32\textheight,width=0.475\textwidth]{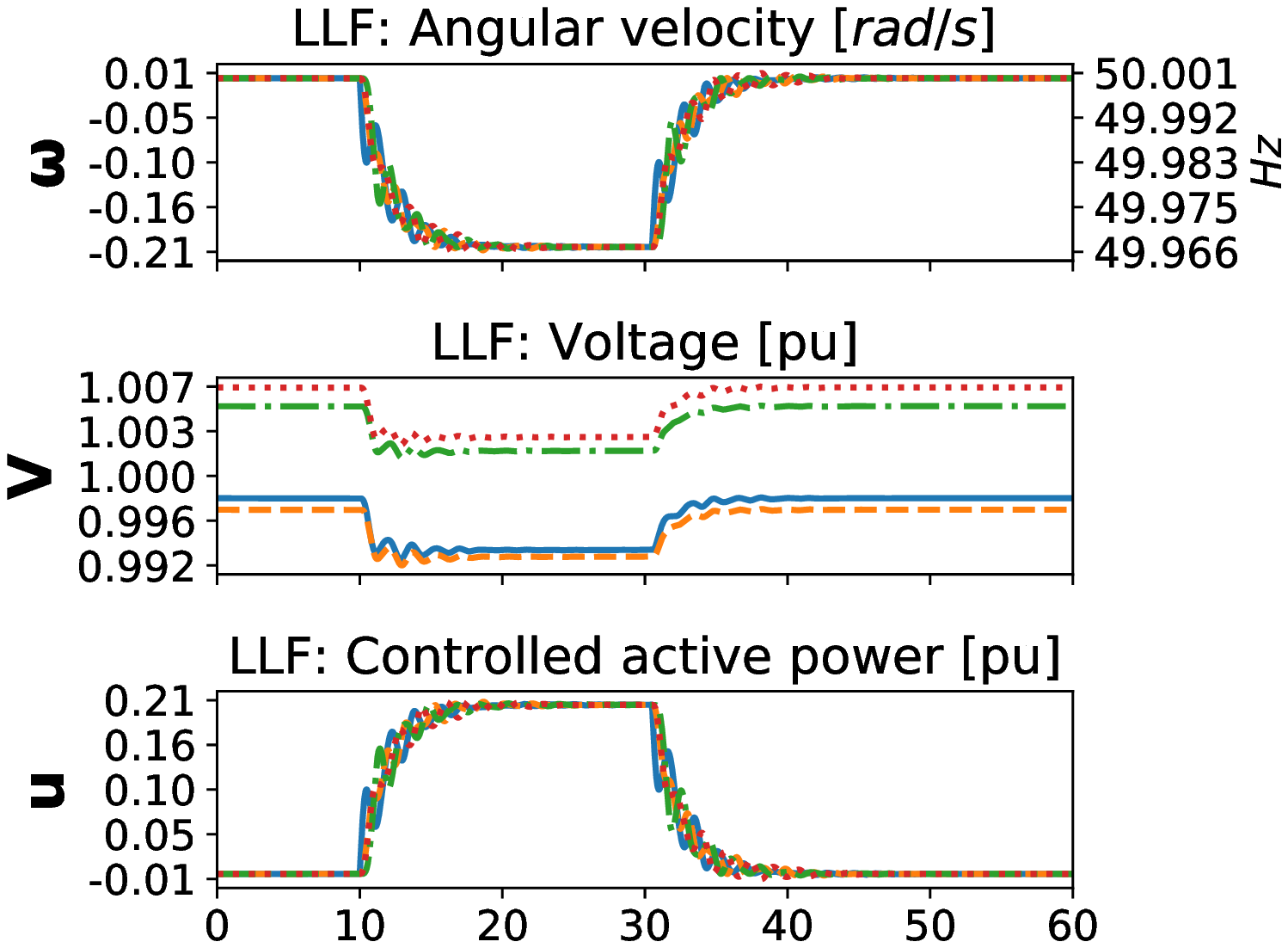}
		\label{fig:LLF_temporary_frequency_voltage_and_control}
	}
	\hfil
	\subfloat[ILF for temporary disturbance]{
		\includegraphics[height=0.32\textheight,width=0.475\textwidth]{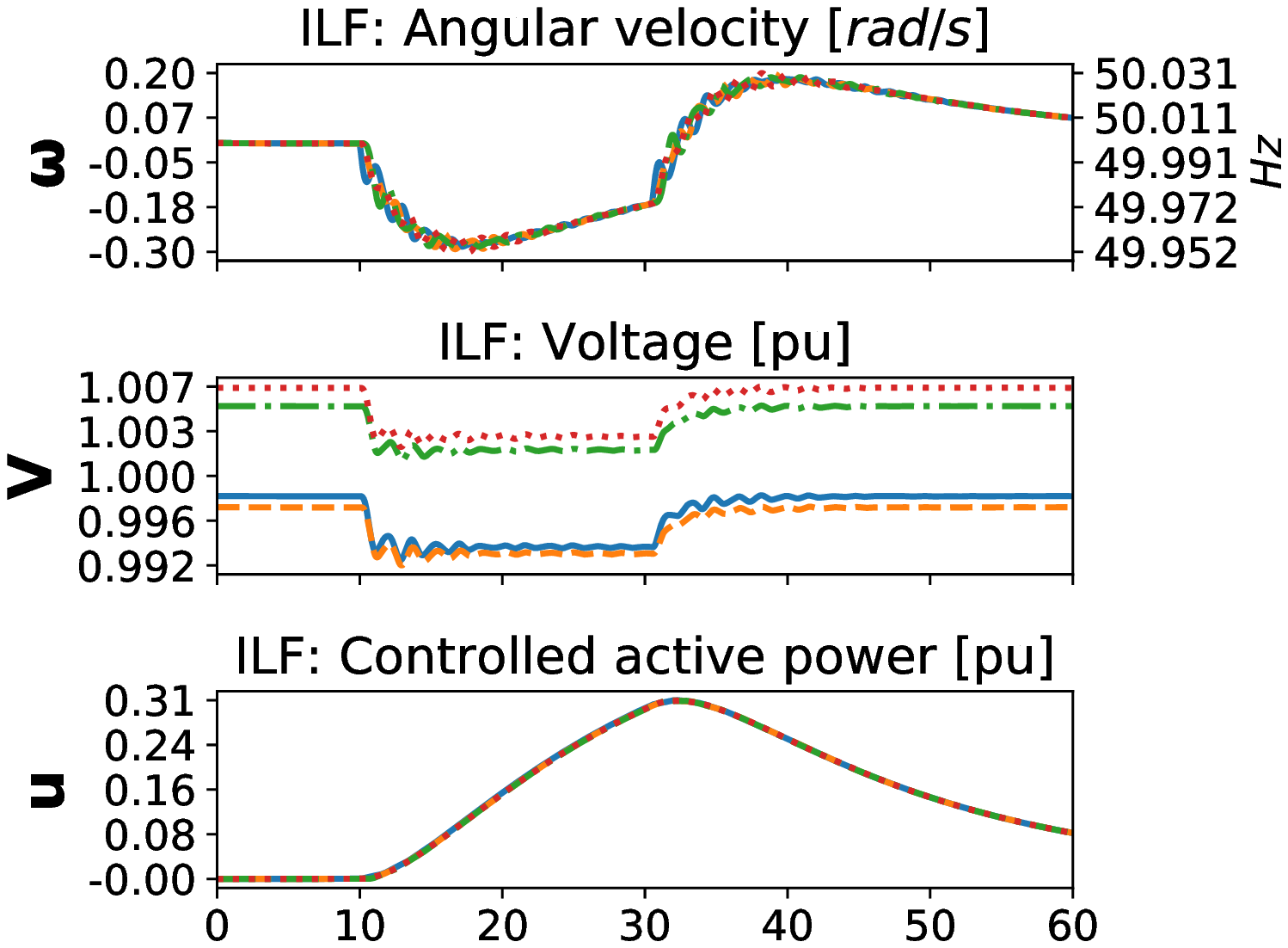}
		\label{fig:ILF_temporary_frequency_voltage_and_control}
	}
	\hfil
	\subfloat[GAB for temporary disturbance]{
		\includegraphics[height=0.32\textheight,width=0.475\textwidth]{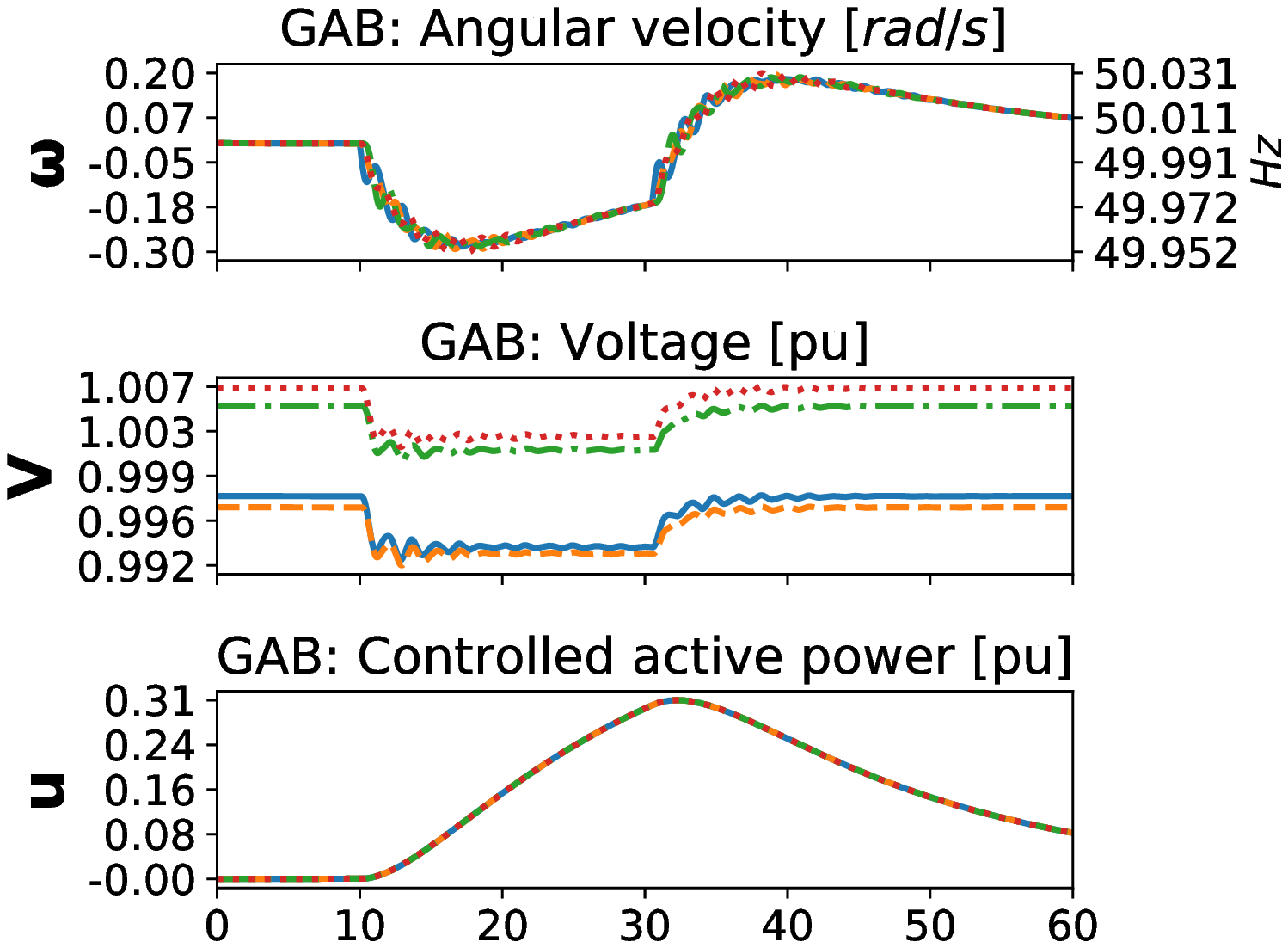}
		\label{fig:GAB_temporary_frequency_voltage_and_control}
	}
	\hfil
	\subfloat[OC for temporary disturbance]{
		\includegraphics[height=0.32\textheight,width=0.475\textwidth]{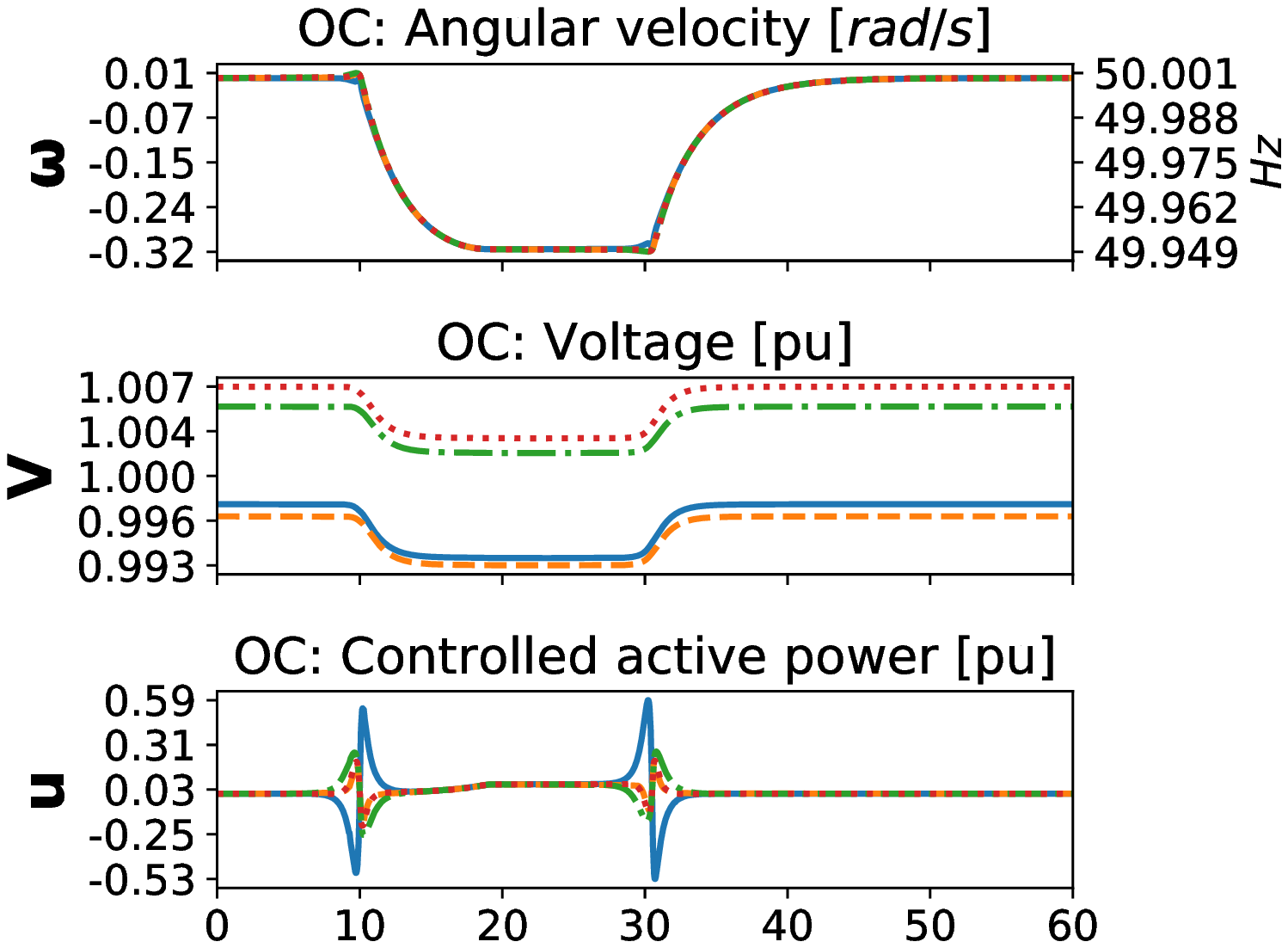}
		\label{fig:OC_temporary_frequency_voltage_and_control}
	}
	\caption{The angular velocity with corresponding frequency values, voltage and controlled power at each node in the test system under the temporary disturbance. Solid, dashed, dash-dotted and dotted lines correspond to nodes 1, 2, 3 and 4 respectively. Each control gradually synchronizes the angular velocities after each change in power by the disturbance. The ILF and GAB controls furthermore try to return the angular velocities to the initial synchronized value. Notice that OC also responds pre-emptively to the disturbance in a significant way.}
	\label{fig:Temporary-Frequency-Voltage-Control}
\end{figure*}

\begin{figure*}[!th]
	\centering
	\subfloat[LLF for temporary disturbance]{
		\includegraphics[height=0.25\textheight,width=0.475\textwidth]{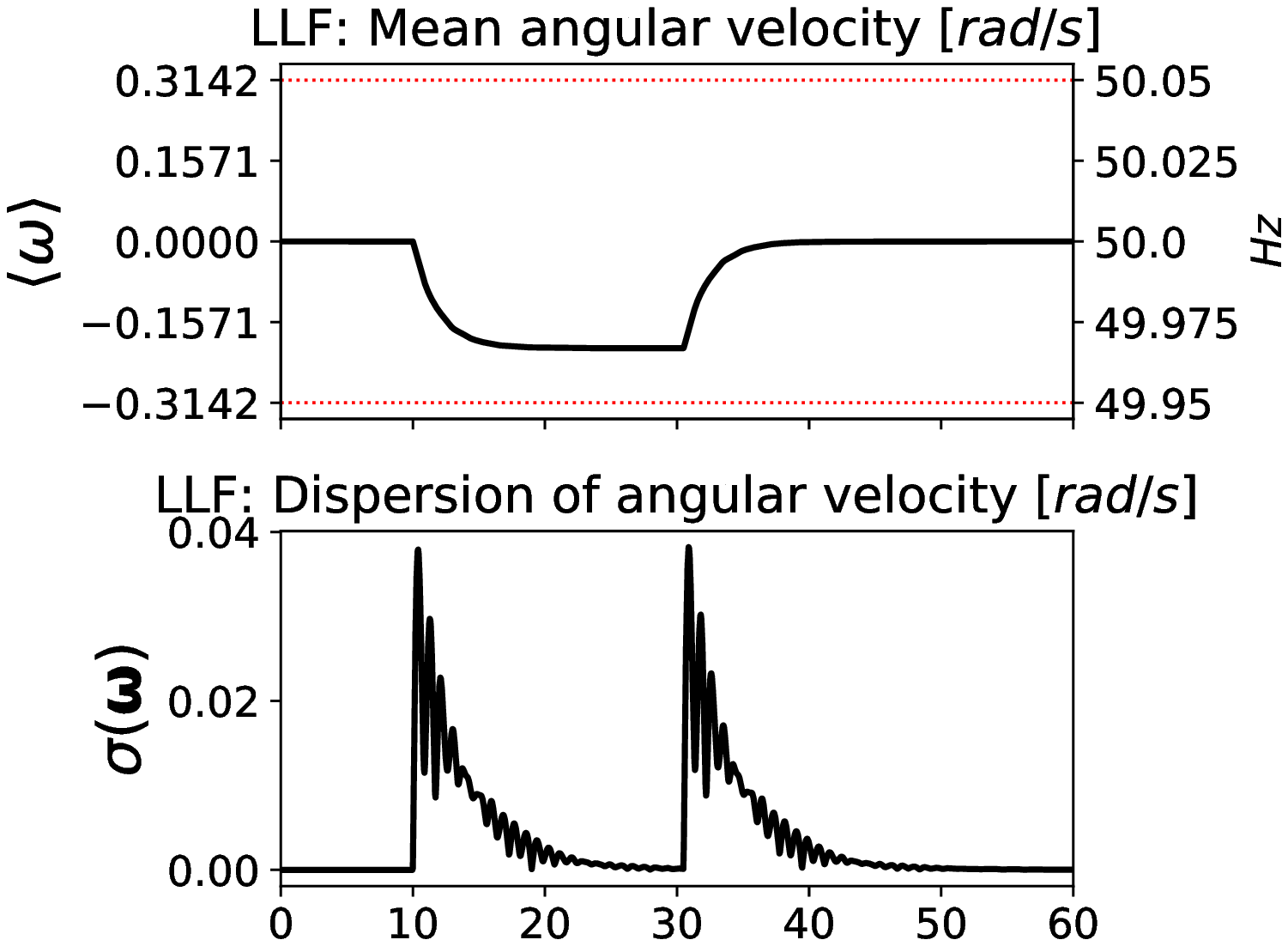}
		\label{fig:LLF_temporary_system_stability}
	}
	\hfil
	\subfloat[ILF for temporary disturbance]{
		\includegraphics[height=0.25\textheight,width=0.475\textwidth]{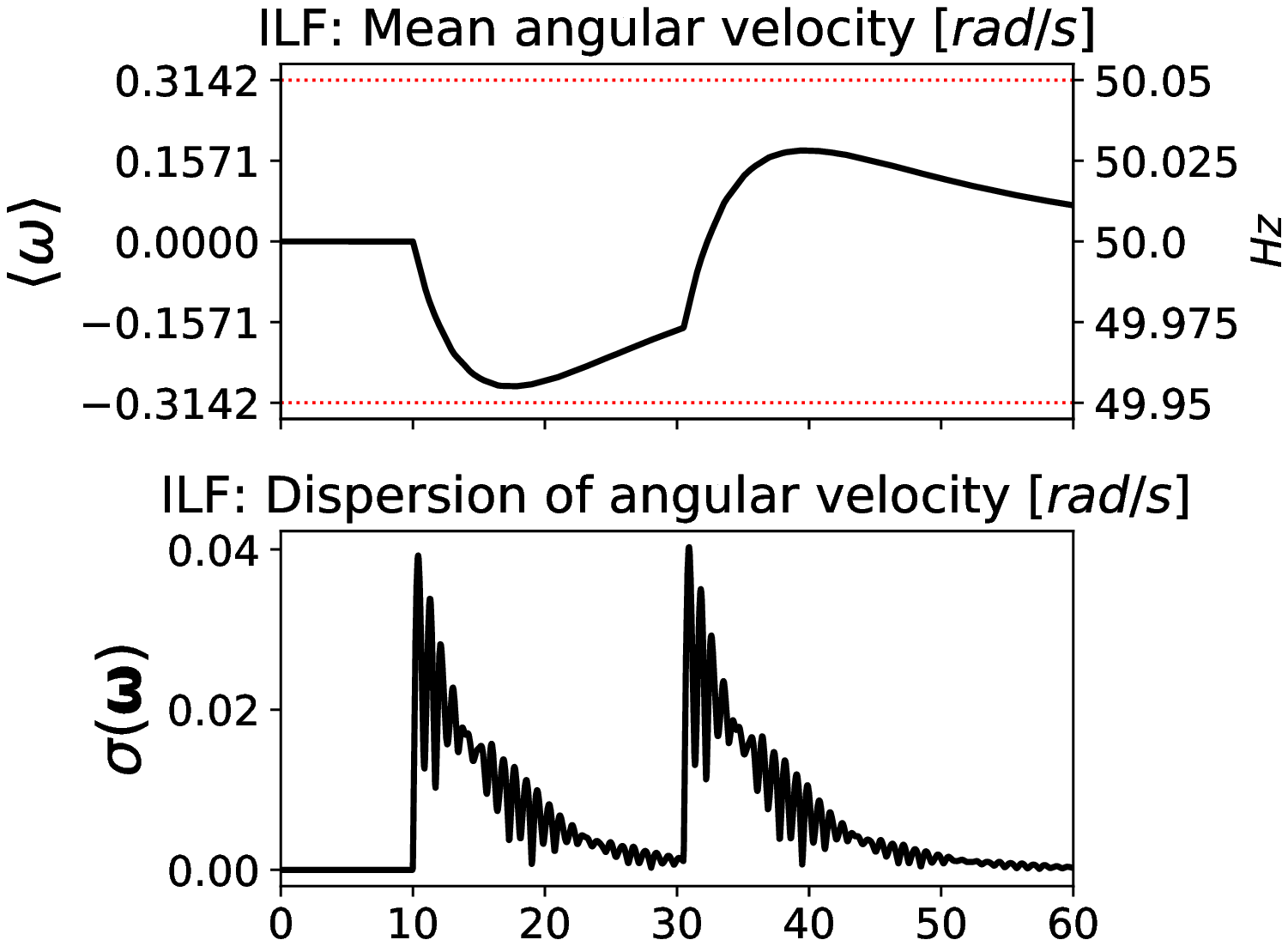}
		\label{fig:ILF_temporary_system_stability}
	}
	\hfil
	\subfloat[GAB for temporary disturbance]{
		\includegraphics[height=0.25\textheight,width=0.475\textwidth]{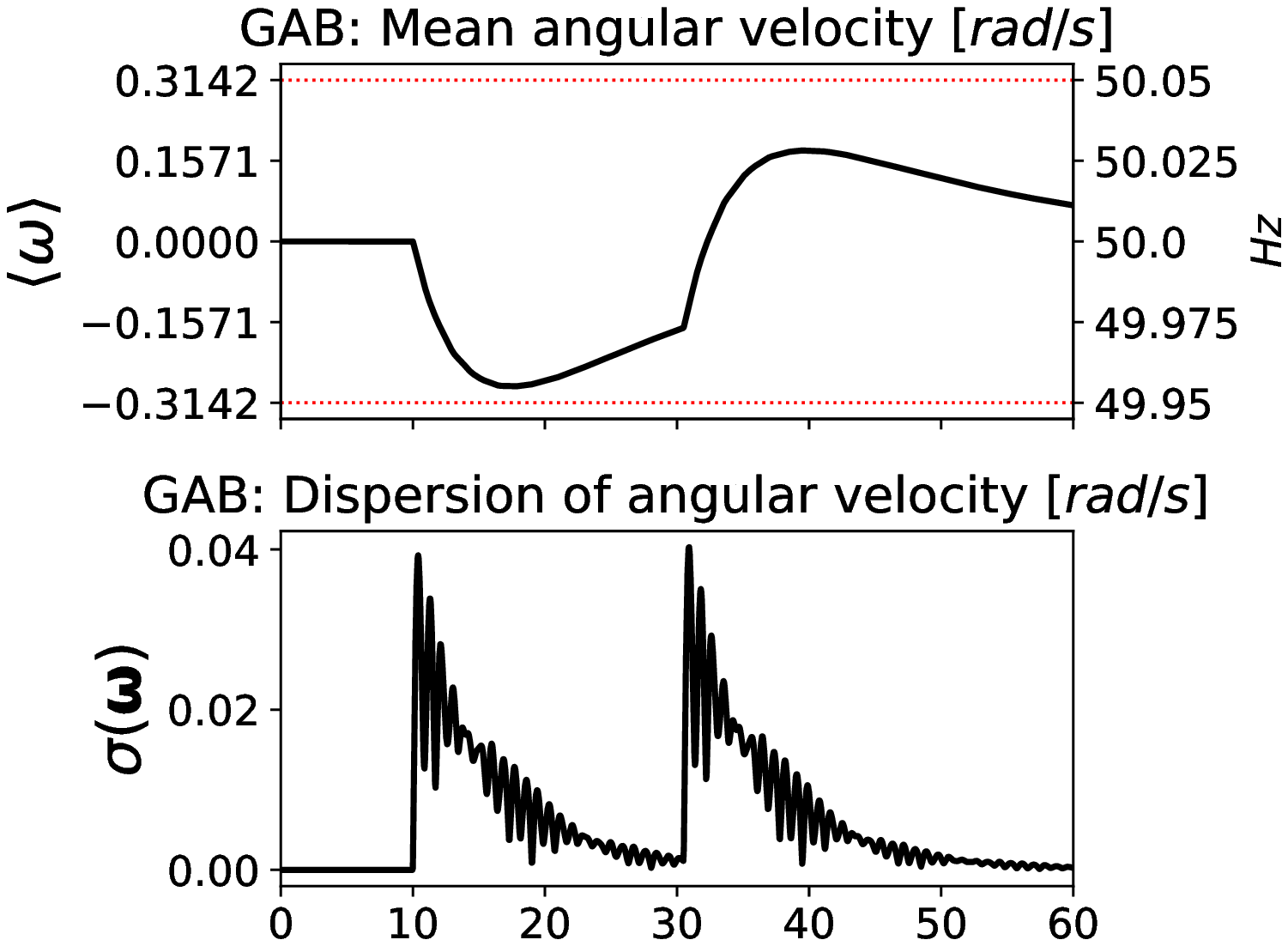}
		\label{fig:GAB_temporary_system_stability}
	}
	\hfil
	\subfloat[OC for temporary disturbance]{
		\includegraphics[height=0.25\textheight,width=0.475\textwidth]{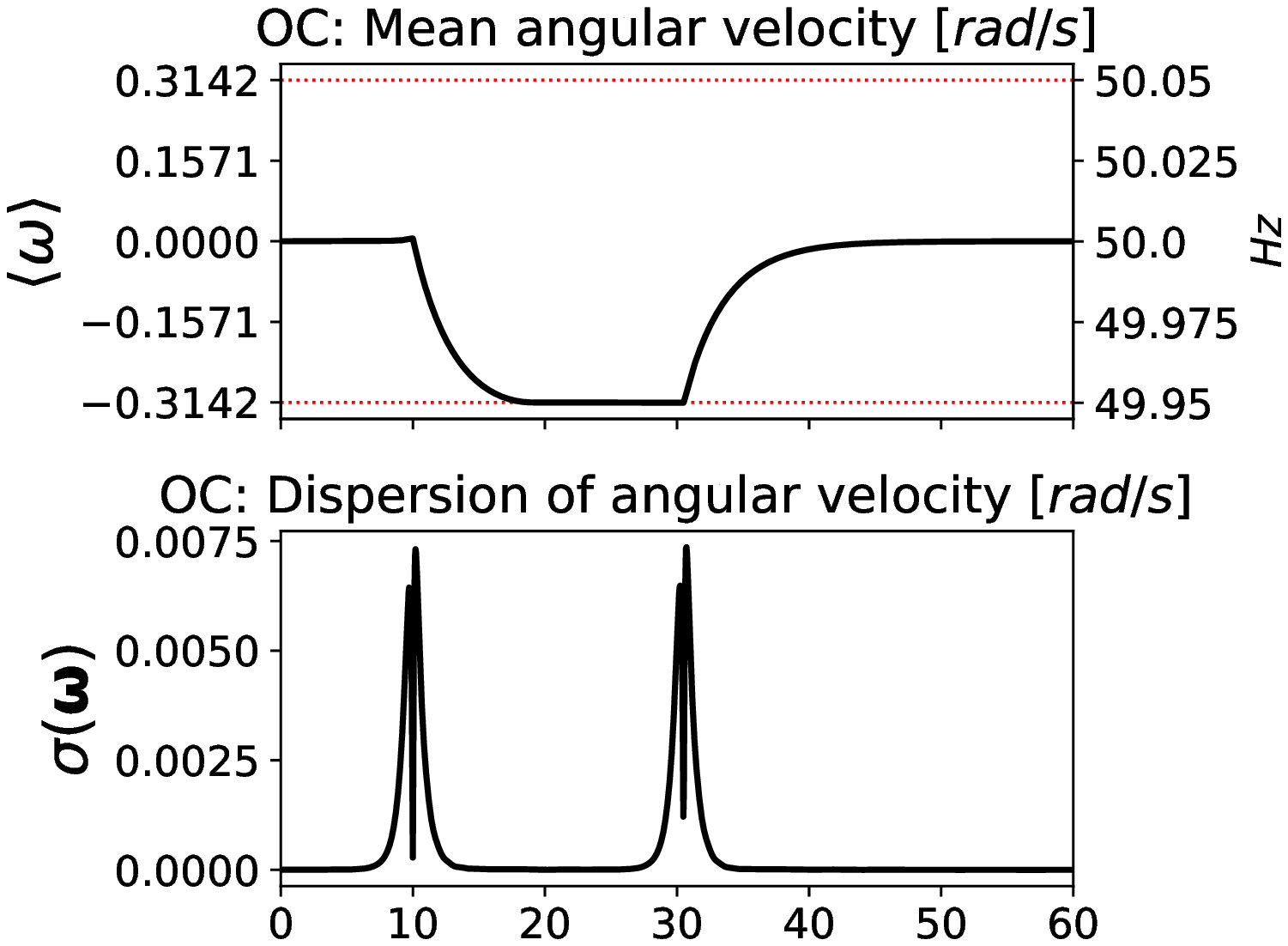}
		\label{fig:OC_temporary_system_stability}
	}
	\caption{Angular velocity mean and deviation in the test system under the temporary disturbance. Red dotted lines show operational limits. Each control keeps the mean angular velocity $\langle \omega \rangle$ within its bounds and gradually synchronizes the system after each change in power by the disturbance. Notice that OC synchronizes the angular velocities to the boundary of its admissible set of values. Furthermore, its pre-emptive responses to the disturbance cause temporary losses of synchronization.}
	\label{fig:Temporary-System-Synchronization}
\end{figure*}

\begin{figure*}[!th]
	\centering
	\includegraphics[width=0.6\textwidth,height=0.3\textheight]{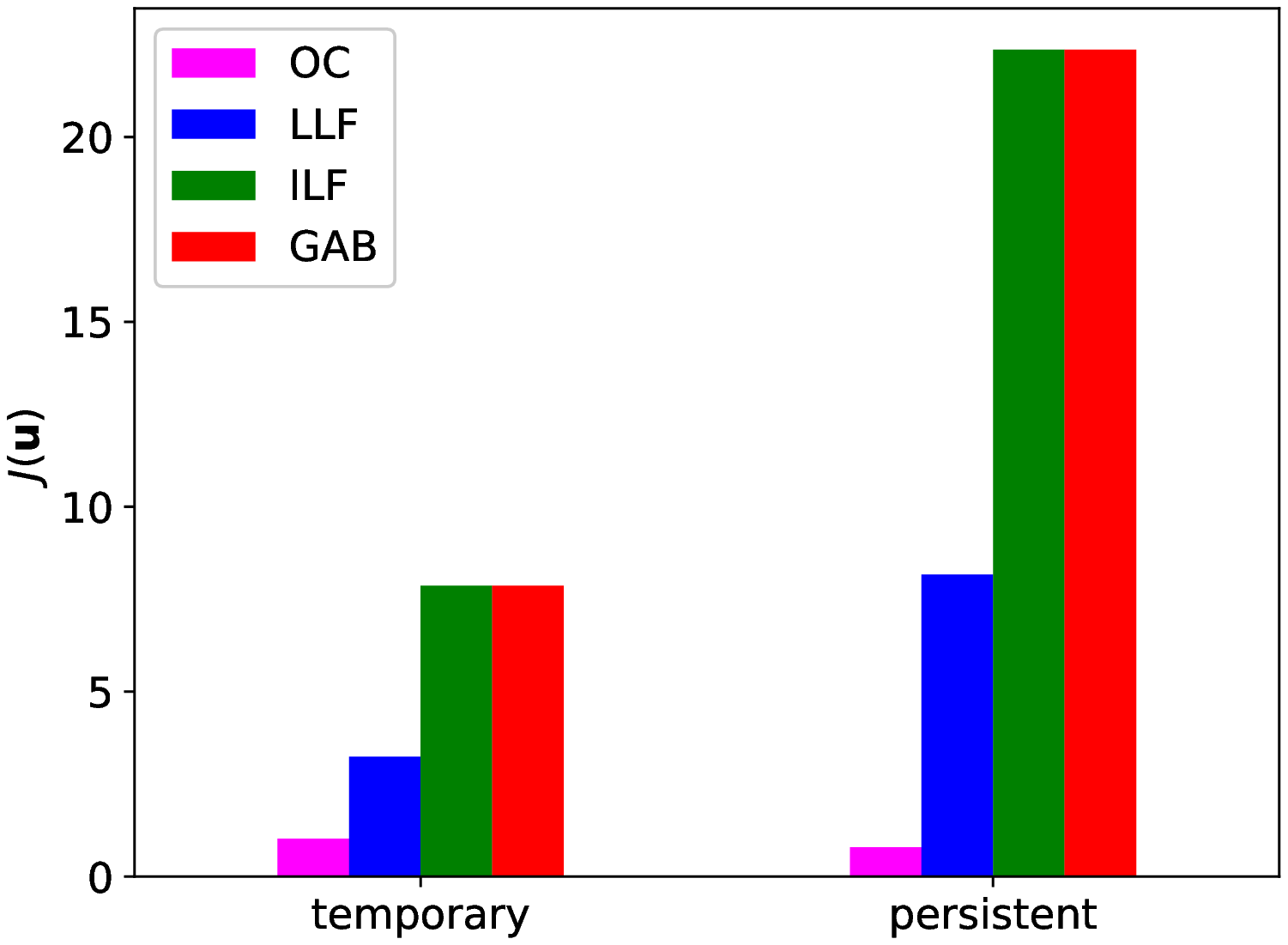}
	\caption{Comparison of control costs for temporary and persistent disturbances (from left to right: LLF, ILF, GAB, OC). The optimal control, OC, keeps the system within operational boundaries at the lowest costs whereas ILF and GAB have the highest costs. We also observe near equal costs for ILF and GAB with other values for the coefficients $\mu$ and $\kappa$ satisfying $\kappa \ge 1$ and $\frac{\mu}{\kappa} = 4$.}
	\label{fig:Control-Costs}
\end{figure*}

Even in the absence of control, the simulated system gradually resynchronizes within the horizon $[0,T]$ with acceptable voltages and, except when the disturbance persists, acceptable angular velocities. We show in Fig.~\ref{fig:Temporary-Frequency-Voltage-Control} trajectories for the controlled active power, angular velocity and voltage under the temporary disturbance, and in Fig.~\ref{fig:Temporary-System-Synchronization} corresponding trajectories for the angular velocity mean and standard deviation. Trajectories under the persistent disturbance display analogous behaviour and are shown in Appendix \ref{Section:Persistent-Trajectories}.

\paragraph*{\bf Linear local frequency (LLF) control}
The LLF control keeps the angular velocities within the given bounds over the control horizon and also synchronizes the system after each change in power by the disturbance. For the persistent disturbance, the angular velocities synchronize near the nadir shown in Fig.~\ref{fig:LLF_temporary_system_stability}. Note that the responses at the nodes become equal as the system synchronizes since the parameters for the control \eqref{eq:Control-Proportional-Frequency} satisfy $\nu_{i} = \nu$ for all $i \in \mathcal{N}$. The displayed control trajectories are oscillatory and dampen while the disturbance $\xi_{1}$ remains constant. However, in separate simulations with larger $\nu$ (not shown) we no longer notice these oscillations. Moreover, when $\nu$ is very large, say $\nu = 100$, the LLF control has a much larger initial response at node 1 that approximates the change in power caused by the disturbance. In this case the angular velocities are also kept much closer to $0$.

\paragraph*{\bf Integral local frequency (ILF) control}
The ILF control also keeps the angular velocities within the bounds over the control horizon and synchronizes the system after each change in power by the disturbance. Moreover, ILF also returns the angular velocities to the initial synchronized value, thereby performing a secondary control action. The displayed control trajectories do not have the oscillations present for the LLF control. However, if $\kappa$ is sufficiently small, then such oscillations can appear, although the angular velocities are kept much closer to $0$.

\paragraph*{\bf Gather-and-broadcast (GAB) control}
The GAB control behaves and performs similarly to ILF as Fig.~\ref{fig:Temporary-System-Synchronization} and results in Table \ref{Table:Control-Coefficient-Values} can attest. In particular, GAB synchronizes the system and performs the secondary control action of returning the frequency to its nominal value.

\paragraph*{\bf Optimal control}
The optimal control causes the mean angular velocity $\langle \omega \rangle$ to follow its natural direction of descent or ascent within the operational limits until a particular level. The angular velocity is then kept at this level whilst the disturbance persists. Additionally, the combined action at the unperturbed nodes is generally of the opposite type to that taken at the perturbed node. That is, when there is an increase (respectively, decrease) in $u_{1}$ there is typically a decrease (respectively, increase) in $\sum_{i=2}^{N}u_{i}$ at the same time. We also notice the following pre-emptive behaviour of the control: shortly before the sudden increase (resp. decrease) in demand at node 1, the optimal control decreases (resp. increases) the active power at this node and simultaneously increases (resp. decreases) the active power at the remaining unperturbed nodes. Consequently, the optimal control uses additional and, in practice, uncertain information about the disturbance in its response that realistic controls may not be able to use. Hence, the optimal controller should always outperform any realistic controller. Finally, we note that the results depend on the parameters selected. For example, if the synchronization loss tolerance is increased from the value $\varepsilon_{1} = 10^{-4}$ (used to generate these results) to $\varepsilon_{1} = 10^{-3}$ we observe oscillations in the control trajectories.

\subsection{Comparison of control costs}

In Fig.~\ref{fig:Control-Costs} we show the cost $J(\boldsymbol{u})$ for the controls LLF, ILF, GAB and OC associated with the trajectories displayed above. While it is clear that OC satisfies the constraints with smallest cost at the lowest sychronization loss (Table \ref{Table:Control-Coefficient-Values}), these costs can depend significantly on the simulation parameters. For example, the LLF cost increases with the coefficient $\nu$ and the OC cost increases as the synchronization loss tolerance $\varepsilon_{1}$ decreases. Notwithstanding this we can explain the disparity between costs for LLF and ILF (or GAB) by the additional secondary control action undertaken by ILF (see Fig.~\ref{fig:Temporary-Frequency-Voltage-Control}). Also, the similarity in costs for the temporary and persistent disturbances corresponding to OC can be attributed to the significant cost of responding pre-emptively to the temporary disturbance in this case.

\section{Conclusion and outlook}\label{Section:Outlook}
In summary, we have introduced and numerically solved an optimal control problem to benchmark different control schemes for power grid transient stability in terms of their economic effectiveness. We investigated three distributed control schemes: linear local frequency \eqref{eq:Control-Proportional-Frequency}, integral local frequency \eqref{eq:Control-Integral-Frequency}, and gather-and-broadcast \eqref{eq:GAB-Control-Global}.

The linear local frequency control acts as a primary response service to keep the grid frequency close to its nominal value. If the control coefficient $\nu_{i}$ in \eqref{eq:Control-Proportional-Frequency} is chosen suitably, for example comparable to the damping parameter at node $i$, then this control can be quite cost effective when compared to the integral frequency and gather-and-broadcast controls. However, we note that the latter controls can also provide secondary response service (see Fig.~\ref{fig:Temporary-Frequency-Voltage-Control}) which the linear local frequency control is not designed for. If the coefficient $\nu_{i}$ for the linear local frequency control is large, this leads to more costly power response profiles that almost exactly counteract the disturbance, at least in the initial response phase. The linear local frequency, integral frequency and gather-and-broadcast controllers can also produce control trajectories with oscillations depending on how their parameters are chosen.

Our results suggest that more efficient controllers distribute the controlled response amongst all nodes in the power grid.  Moreover, this response need not be homogeneous throughout the network, but could simultaneously involve incremental actions (net increase in power) at some nodes and decremental ones (net decrease in power) at others. Trajectories associated with the optimal control show that as it changes the net active power, the mean angular velocity follows its natural direction of descent, or ascent as appropriate, within the operational limits until a point is reached, possibly at the boundary, at which the power grid is synchronized and active power is balanced within the network.

A response like the one exhibited by the optimal control apparently requires additional information about the disturbance that is likely to be uncertain. Nevertheless, for events that are planned or will occur with very high probability at an anticipated future time, information about the disturbance can be incorporated in the control system's initial response, and a simple distributed or decentralized control such as those we investigated can be used to smooth out additional unknown perturbations. Designing {\it optimal} distributed controllers is the subject of ongoing work (see \cite{Molzahn2017,Stegink2017}, for instance) and decentralized stochastic control (see \cite{Mahajan2016,Charalambous2017,Singh2016}), which generalizes our methodology by incorporating uncertainties and different information structures amongst multiple controllers, is likely to become an important theoretical tool for understanding how these controllers work. Finally, while the numerical results presented here were obtained for a specific four-node network, they provide useful heuristics for more realistic and larger networks. Overall, our results contribute insight into the process of designing and benchmarking secure and cost-efficient controllers for the power system.



\newpage
\appendices
\section{Tables of parameter values}\label{Section:Parameter-Tables}

\begin{table}[!h]
	\renewcommand{\arraystretch}{1.3}
	\caption{Steady state values and parameters for the power grid model \eqref{eq:3rd-order-controlled-swing-equation} used in the simulations, based on \cite[p.~251]{Trip2016}. The net injection $P_{in,i}$ is obtained from the mechanical power $P_{m,i}$ and aggregate load $P_{l,i}$ by $P_{in,i} = P_{m,i} - P_{l,i}$. Line susceptance values $B_{i,j}$ other than those listed are equal to $0$ except $B_{1,2} = B_{2,1} = 34.13$, $B_{1,4} = B_{4,1} = 28$, $B_{2,3} = B_{3,2} = 44.1$ and $B_{3,4} = B_{4,3} = 22.1$.}
	\label{Table:Simulation-Parameter-Values}
	\centering
	\begin{tabular}{ccccc}
		Parameter [units] & Node 1 & Node 2 & Node 3 & Node 4 \\\hline
		$M_{i}$ $[s^{2}]$ & $5.22$ & $3.98$ & $4.49$ & $4.22$ \\
		$D_{i}$ [pu] & $1.60$ & $1.22$ & $1.38$ & $1.42$ \\
		$E_{f,i}$ [pu] & $7.01$ & $6.09$ & $6.29$ & $6.67$ \\
		$T'_{do,i}$ $[s]$ & $5.54$ & $7.41$ & $6.11$ & $6.22$ \\
		$X_{d,i}$ [pu] & $1.84$ & $1.62$ & $1.80$ & $1.94$ \\
		$X'_{d,i}$ [pu] & $0.25$ & $0.17$ & $0.36$ & $0.44$ \\
		$B_{i,i}$ [pu] & $-66.1$ & $-82.2$ & $-69.6$ & $-53.6$ \\
		$P_{m,i}$ [pu] & $1.1$ & $1.4$ & $0.8$ & $2.2$ \\
		$P_{l,i}$ [pu] & $2.0$ & $1.0$ & $1.5$  & $1.0$ \\
		$P_{in,i}$ [pu] & $-0.9$ & $0.4$ & $-0.7$ & $1.2$ \\
		$\bar{\theta}_{i}$ [$rad$] & $0.0911$ & $0.0973$ & $0.0930$ & $0.115$ \\
		$\bar{\omega}_{i}$ [$rad \cdot s^{-1}$] & $0$ & $0$ & $0$ & $0$ \\
		$\bar{V}_{i}$ [pu] & $0.998$ & $0.997$ & $1$ & $1$ \\\hline
	\end{tabular}
\end{table}

\begin{table}[!h]
	\renewcommand{\arraystretch}{1.3}
	\caption{Parameter values for the control problem. Vectors and matrices are denoted in boldface. The tolerance $\varepsilon_{1}$ for the synchronization constraint is set larger than the other tolerances to allow for the loss of synchronization around the occurrence of a disturbance.}
	\centering
	\begin{tabular}{ccc}
		Parameter & Value & Units \\\hline
		$T$ & $60$ & $s$ \\
		$\omega_{min}$, $\omega_{max}$ & $-\frac{\pi}{10}$, $\frac{\pi}{10}$ & $rad \cdot s^{-1}$ \\
		$\boldsymbol{\lambda}$ & $\boldsymbol{1}$ & 1 \\
		$\varepsilon_{1}$ & $10^{-4}$ & 1 \\
		$\varepsilon_{2}$,\ldots,$\varepsilon_{6}$ & $10^{-10}$  & 1 \\
		$\boldsymbol{U^{min}}$, $\boldsymbol{U^{max}}$ & $\boldsymbol{-5}$, $\boldsymbol{5}$ & pu \\
		$\boldsymbol{V^{min}}$, $\boldsymbol{V^{max}}$ & $\boldsymbol{0.94}$, $\boldsymbol{1.06}$ & pu \\
		\hline
	\end{tabular}
\end{table}

\section{The control parametrization method}\label{Appendix:Control-Parametrization-Method}
The following description of the control parametrization method is summarized from the textbook \cite{Teo1991}. Further extensions to this method can be found in the survey \cite{Lin2013}. Let $S^{p}$, where $p \ge 1$ is an integer, denote a finite subset of the control horizon $[0,T]$ consisting of $n_{p} + 1$ {\it partitioning points} $t^{p}_{0},\ldots,t^{p}_{n_{p}}$,
\[
t^{p}_{0} = 0,\; t^{p}_{n_{p}} = T,\;\text{and}\; t^{p}_{k-1} < t^{p}_{k} \;\; \text{for}\; k = 1,\ldots,n_{p}.
\]
An increasing sequence of sets $\{S^{p}\}_{p = 1}^{\infty}$ is formed by taking successive refinements of partitioning points, and these sets should become dense in $[0,T]$ as $p$ tends to infinity,
\[
\lim_{p \to \infty}\max_{k = 1,\ldots,n_{p}}|t^{p}_{k} - t^{p}_{k-1}| = 0.
\]
For instance, we can use equidistant partitioning points, $t^{p}_{k} = \frac{k}{n_{p}}T$ for $k = 0,\ldots,n_{p}$, with the ratio $\frac{n_{p+1}}{n_{p}}$, $p \ge 1$, being a constant integer that is greater than 1 (a common choice is $\frac{n_{p+1}}{n_{p}} = 2$). We define $\mathcal{U}^{p}$ as the subset of control variables $\mathbf{u}^{p} \in \mathcal{U}$ that are piecewise constant and consistent with $S^{p}$ in the following sense,
\[
u^{p}_{i}(t) = \sum_{k = 1}^{n_{p}}u_{i}^{p,k}\mathbf{1}_{[t^{p}_{k-1},t^{p}_{k})}(t),\quad u_{i}^{p,k} \in U_{i},\; i \in \mathcal{N}.
\]
Each control $\mathbf{u}^{p}$ is parametrized by an element $\boldsymbol{U}^{p}$ of the $(N \times n_{p})$-dimensional space $\mathbb{U}^{p} = \prod_{k = 1}^{n_{p}}\left(\prod_{i = 1}^{N}U_{i}\right)$, where $\boldsymbol{U}^{p} = \{\boldsymbol{u}^{p}_{k}\}_{k = 1}^{n_{p}}$ and $\boldsymbol{u}^{p}_{k} = (u_{1}^{p,k},\ldots,u_{N}^{p,k})$, This induces equivalent state dynamics $\tilde{\boldsymbol{f}}$, costs $\tilde{J}$ and constraints $\tilde{C}_{\eta}$ that are dependent on the parameter $\boldsymbol{U}^{p}$,
\begin{align*}
\dot{\boldsymbol{x}}(t) & = \tilde{\boldsymbol{f}}(t,\boldsymbol{x}(t),\boldsymbol{U}^{p}) = \boldsymbol{f}(t,\boldsymbol{x}(t),\boldsymbol{u}^{p}(t)),\\ \tilde{J}(\boldsymbol{U}^{p}) & = J(\boldsymbol{u}^{p}), \\
\tilde{C}_{\eta}(\boldsymbol{U}^{p}) & = C_{\eta}(\boldsymbol{u}^{p}).
\end{align*}
An approximate solution to the {\it infinite dimensional} optimal control problem \eqref{Problem:Optimal-Control} is obtained by solving the following non-linear {\it finite dimensional} optimization problem.

\begin{quotation}
	\noindent{\bf Problem:}
	\begin{gather*}
	\text{minimize}\; \tilde{J}(\boldsymbol{U}^{p})\; \text{subject to:} \nonumber \\
	\begin{split}
	i) & \quad
	\dot{\boldsymbol{x}}(t)  = \tilde{\boldsymbol{f}}(t,\boldsymbol{x}(t),\boldsymbol{U}^{p}), \quad
	\boldsymbol{x}(0) = \boldsymbol{x}_{0}; \\
	ii) & \quad \boldsymbol{U}^{p} \in \mathbb{U}^{p}; \\
	iii) & \quad \tilde{C}_{\eta}(\boldsymbol{U}^{p}) \le \varepsilon_{\eta}  \enskip \text{for} \enskip \eta = 1,\ldots,N+2.
	\end{split}
	\end{gather*}
\end{quotation}
An optimization algorithm such as sequential quadratic programming can be used to solve this approximate problem. Such optimization algorithms are typically iterative, and the main computations carried out during each iteration are outlined below (see Section 6.6 of \cite{Teo1991} for further details and \cite{Martyr2018} for an implementation):
\begin{enumerate}
	\item Obtain a trajectory for the state variable $\boldsymbol{x}$ by numerically integrating its dynamics forward in time on the partitioning points $S^{p}$.
	\item Evaluate the cost $\tilde{J}(\boldsymbol{U}^{p})$ and constraints $\tilde{C}_{\eta}(\boldsymbol{U}^{p})$ using numerical integration.
	\item Compute the gradients of the cost $\tilde{J}(\boldsymbol{U}^{p})$ and constraints $\tilde{C}_{\eta}(\boldsymbol{U}^{p})$ according to the formulas given in Section 6.6 of \cite{Teo1991}.
\end{enumerate}
The gradient of the cost $\tilde{J}(\boldsymbol{U}^{p})$, for example, involves computation of the gradient of a {\it Hamiltonian function} $\tilde{\mathcal{H}}$ with respect to the parameter $\boldsymbol{U}^{p}$,
\[
\frac{\partial\tilde{J}(\boldsymbol{U}^{p})}{\partial \boldsymbol{U}^{p}} = \int_{0}^{T}\frac{\partial \tilde{\mathcal{H}}(t,\boldsymbol{x}(t), \boldsymbol{U}^{p},\boldsymbol{z}(t))}{\partial \boldsymbol{U}^{p}}{d}t,
\]
where $\boldsymbol{z}$ is the {\it costate variable} associated to the cost. The Hamiltonian is defined by,
\[
\begin{split}
\tilde{\mathcal{H}}(t,\boldsymbol{x}(t),\boldsymbol{U}^{p},\boldsymbol{z}(t)) = {} & L(t,\boldsymbol{x}(t),\boldsymbol{u}^{p}(t)) \\
& + 
\boldsymbol{z}(t) \cdot \boldsymbol{f}(t,\boldsymbol{x}(t),\boldsymbol{u}^{p}(t)),
\end{split}
\]
where $L$ is the cost rate function in \eqref{eq:Optimal-Control-Lagrange-Criterion} and $\cdot$ is the dot product. Dynamics for this costate variable are given by,
\[
\begin{cases}
\dot{\boldsymbol{z}}(t) = {\displaystyle -\frac{\partial \tilde{\mathcal{H}}(t,\boldsymbol{x}(t),\boldsymbol{U}^{p},\boldsymbol{z}(t))}{\partial \boldsymbol{x}}}\\
\boldsymbol{z}(T) = \boldsymbol{0},
\end{cases}
\]
and this differential equation is solved numerically {\it backwards in time} given a trajectory for $\boldsymbol{x}$. Costate variables for the constraints are defined similarly, but their boundary values at $T$ are non-zero in general due to the presence of terminal costs.
\onecolumn
\section{Simulations under the persistent disturbance}\label{Section:Persistent-Trajectories}

\begin{figure*}[!h]
	\centering
	\subfloat[LLF for persistent disturbance]{
		\includegraphics[height=0.32\textheight,width=0.475\textwidth]{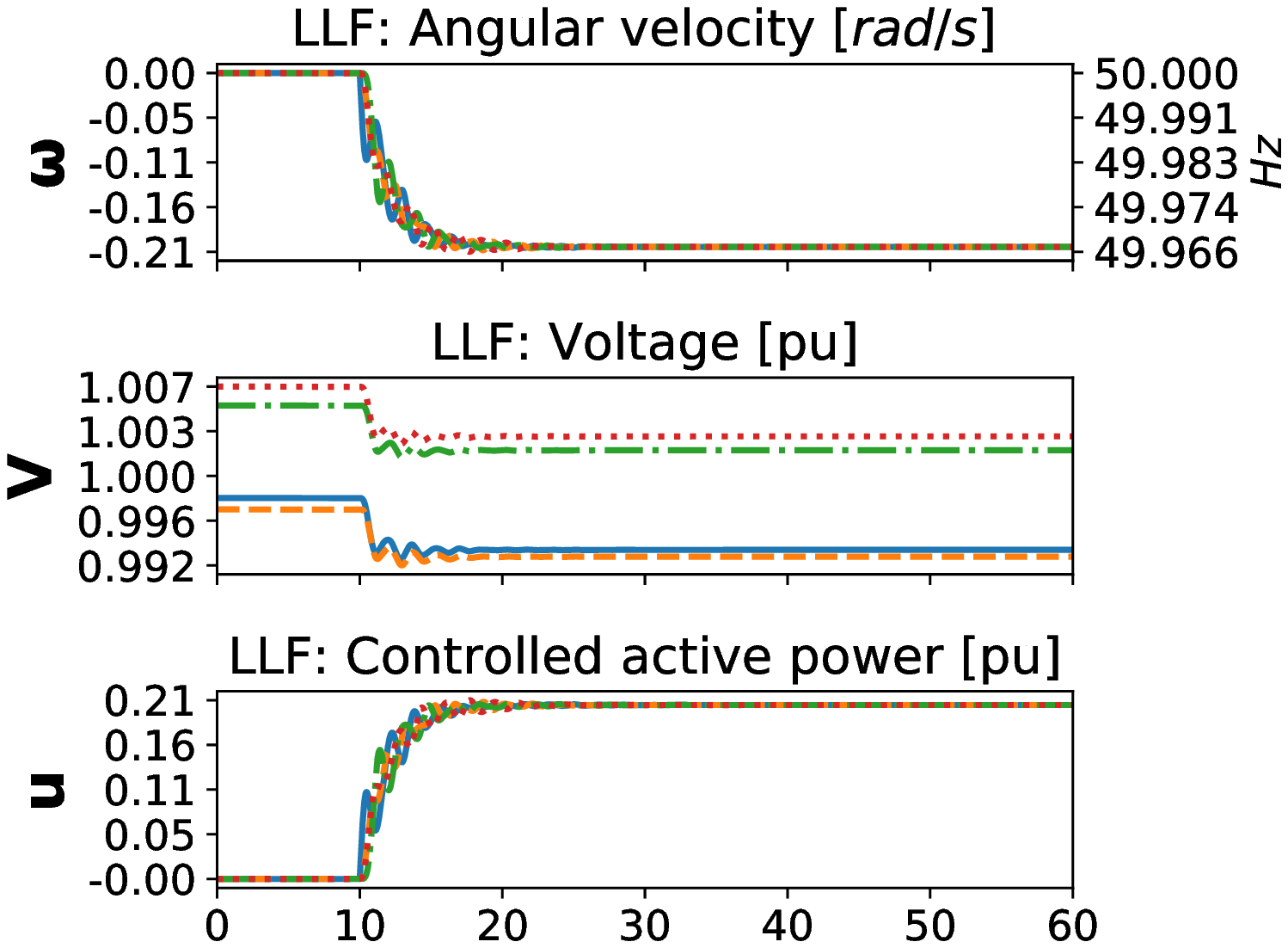}
		\label{fig:LLF_persistent_frequency_voltage_and_control}
	}
	\hfil
	\subfloat[ILF for persistent disturbance]{
		\includegraphics[height=0.32\textheight,width=0.475\textwidth]{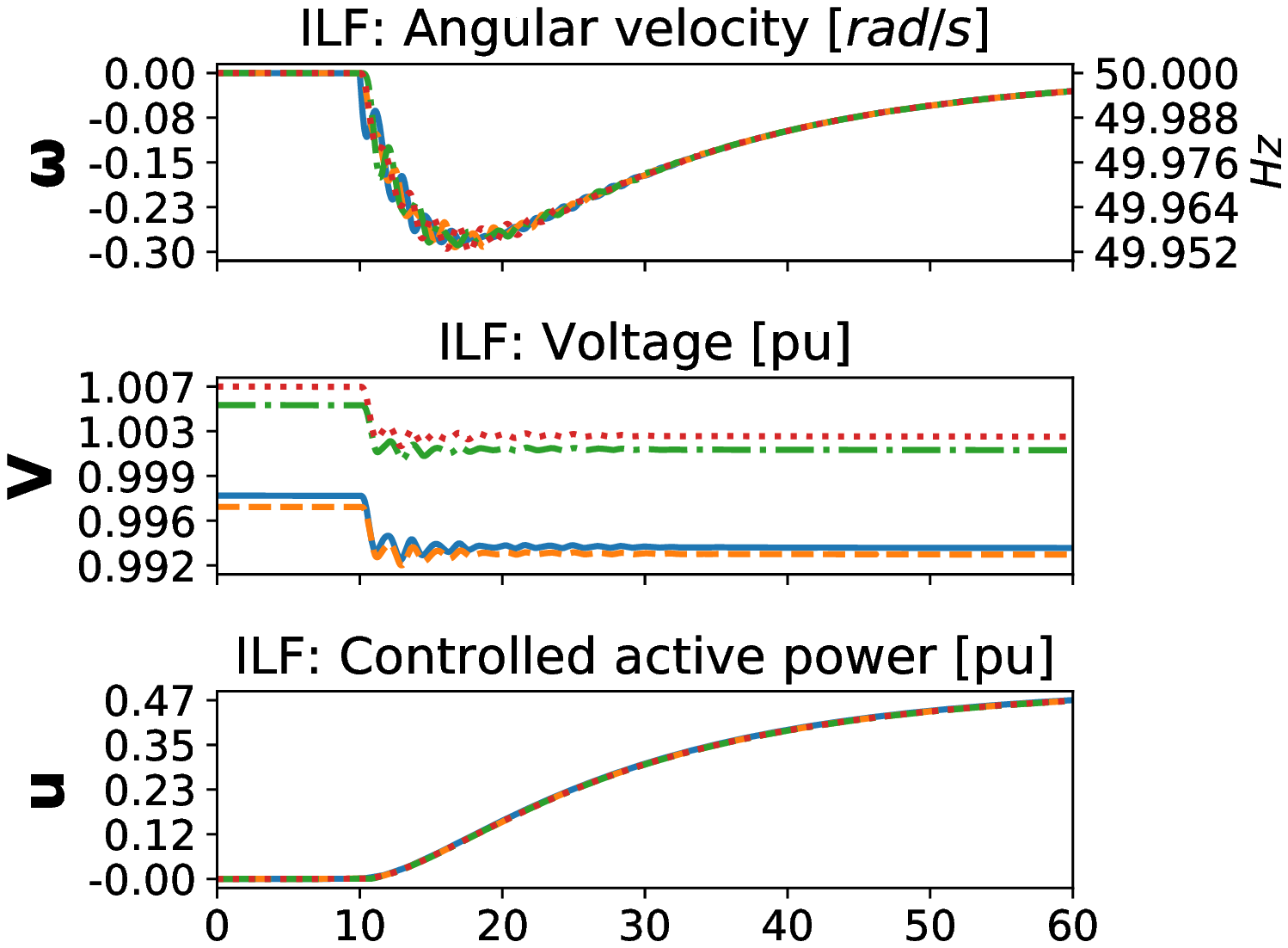}
		\label{fig:ILF_persistent_frequency_voltage_and_control}
	}
	\hfil
	\subfloat[GAB for persistent disturbance]{
		\includegraphics[height=0.32\textheight,width=0.475\textwidth]{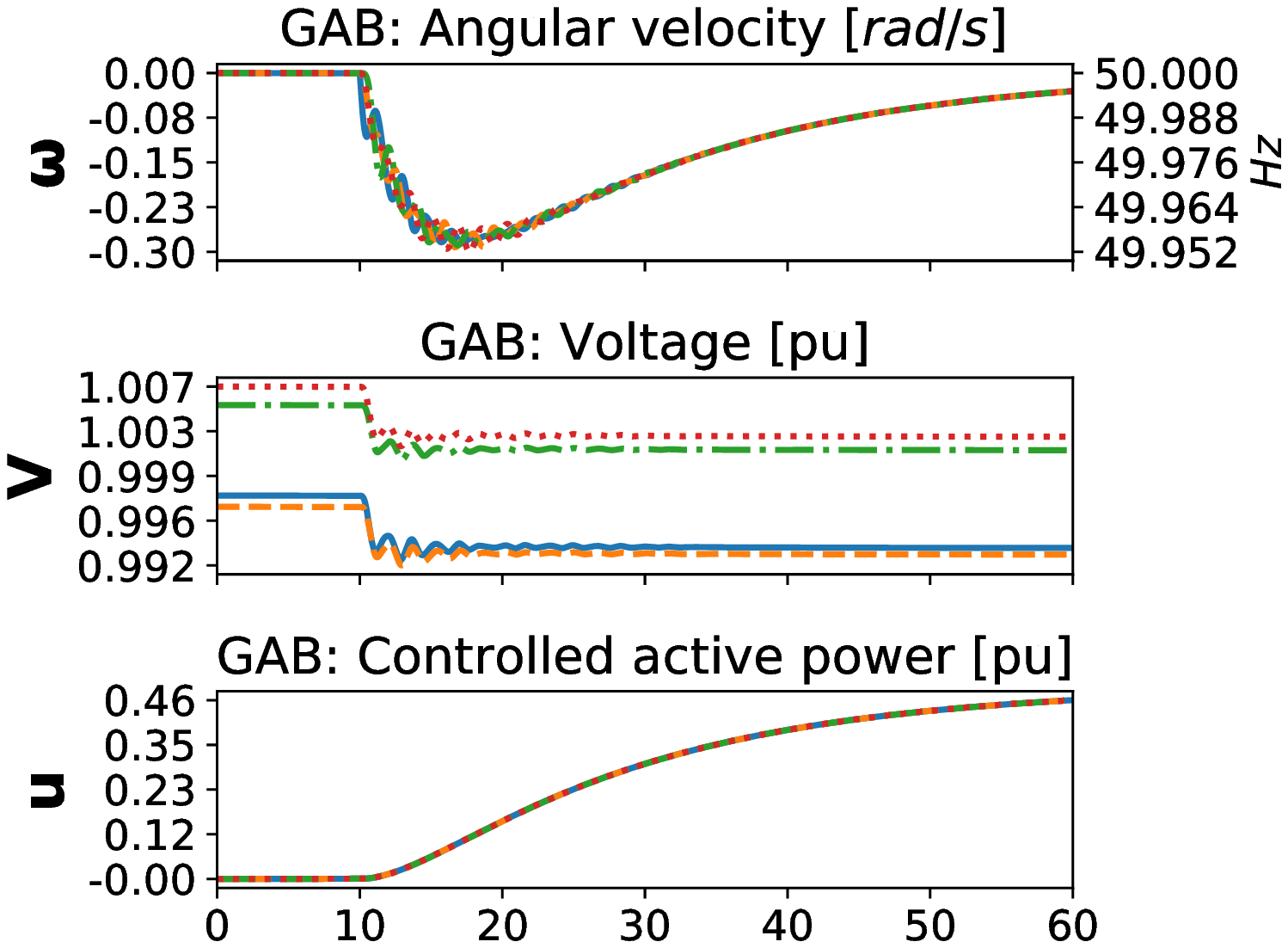}
		\label{fig:GAB_persistent_frequency_voltage_and_control}
	}
	\hfil
	\subfloat[OC for persistent disturbance]{
		\includegraphics[height=0.32\textheight,width=0.475\textwidth]{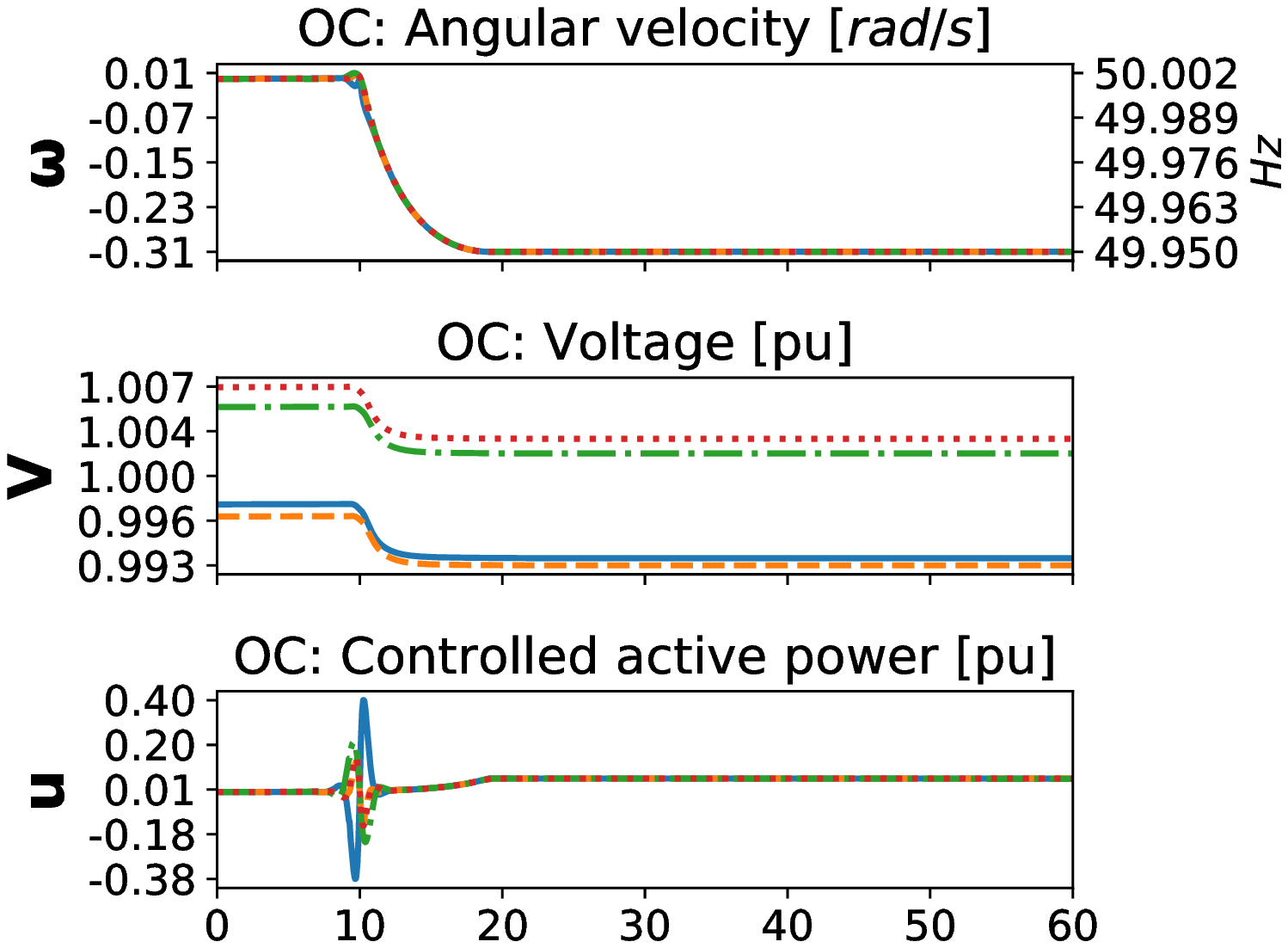}
		\label{fig:OC_persistent_frequency_voltage_and_control}
	}
	\caption{The angular velocity with corresponding frequency values, voltage and controlled power at each node in the test system under the persistent disturbance. Solid, dashed, dash-dotted and dotted lines correspond to nodes 1, 2, 3 and 4 respectively. Each control gradually synchronizes the angular velocities after each change in power by the disturbance. The ILF and GAB controls furthermore try to return the angular velocities to the initial synchronized value.}
	\label{fig:Persistent-Frequency-Voltage-Control}
\end{figure*}

\begin{figure*}[!h]
	\centering
	\subfloat[LLF for persistent disturbance]{
		\includegraphics[height=0.25\textheight,width=0.475\textwidth]{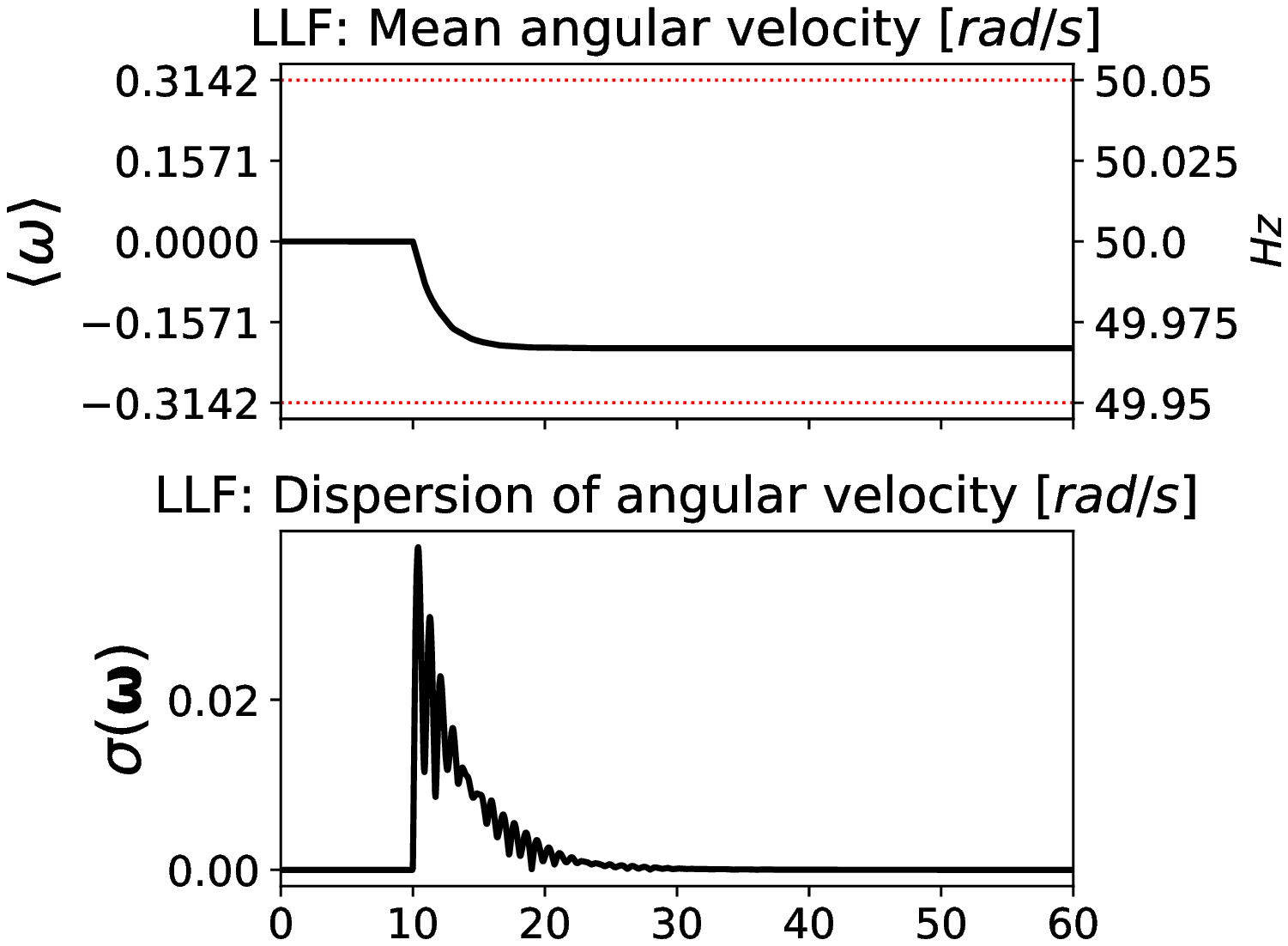}
		\label{fig:LLF_persistent_system_stability}
	}
	\hfil
	\subfloat[ILF for persistent disturbance]{
		\includegraphics[height=0.25\textheight,width=0.475\textwidth]{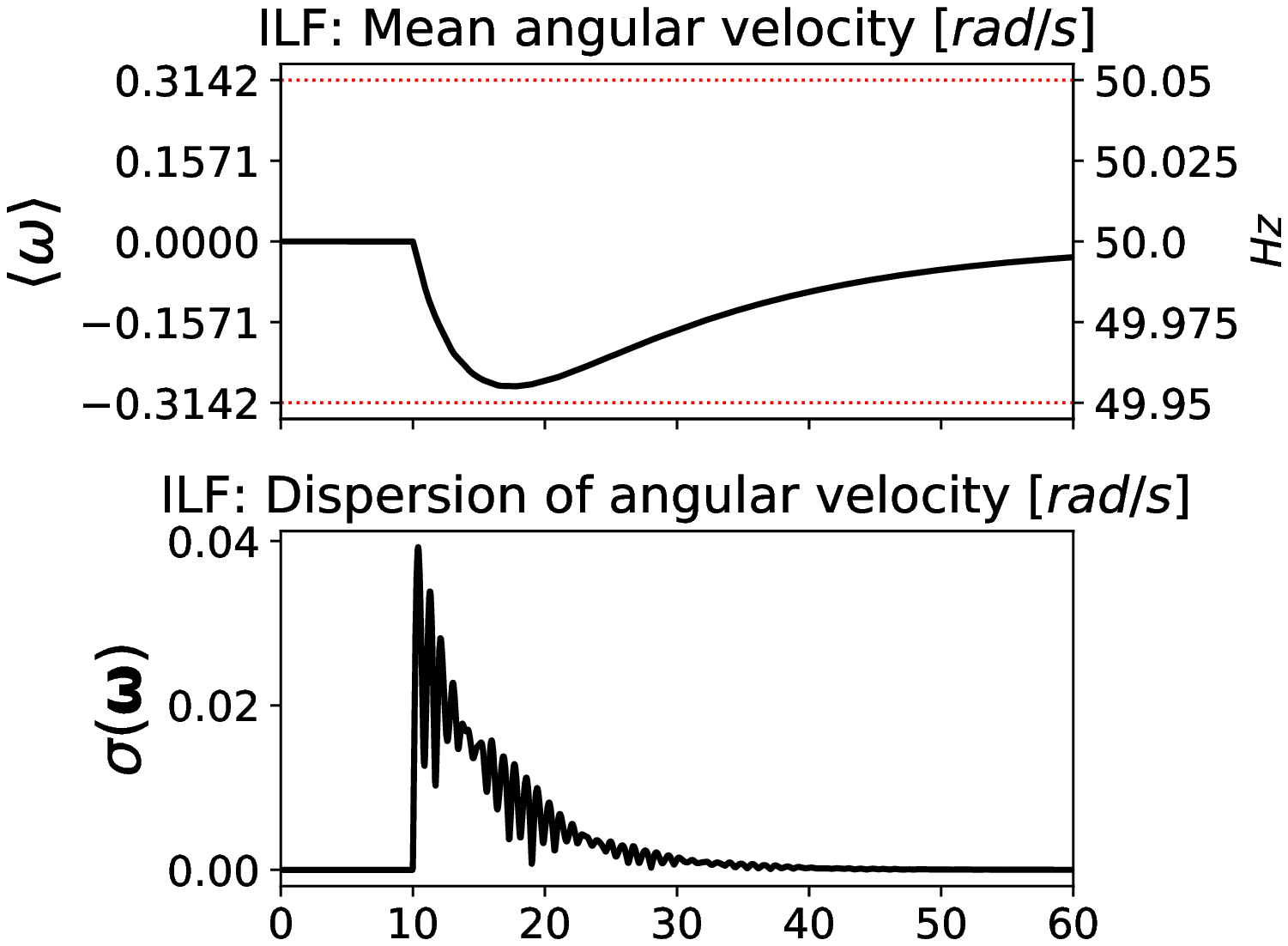}
		\label{fig:ILF_persistent_system_stability}
	}
	\hfil
	\subfloat[GAB for persistent disturbance]{
		\includegraphics[height=0.25\textheight,width=0.475\textwidth]{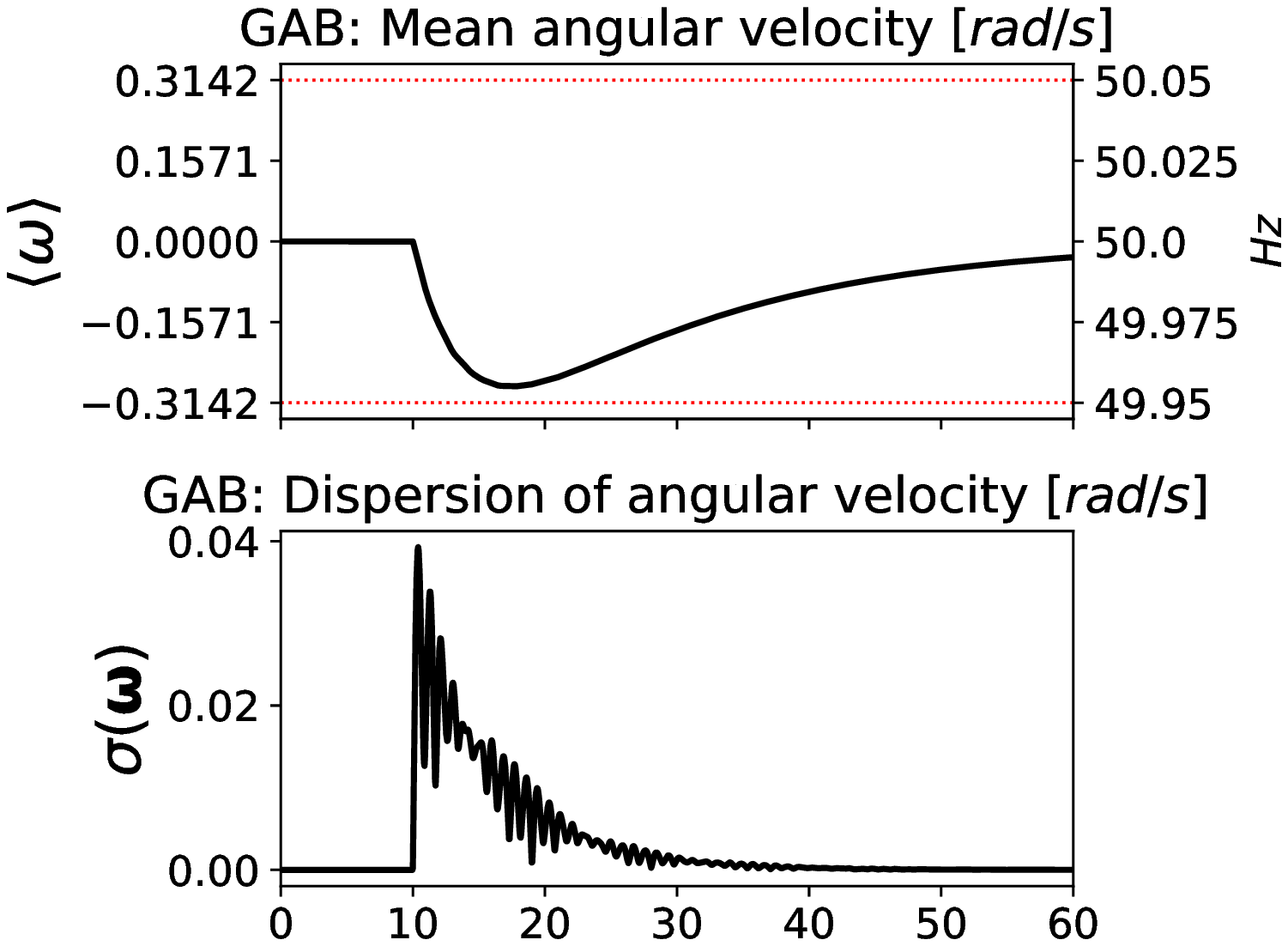}
		\label{fig:GAB_persistent_system_stability}
	}
	\hfil
	\subfloat[OC for persistent disturbance]{
		\includegraphics[height=0.25\textheight,width=0.475\textwidth]{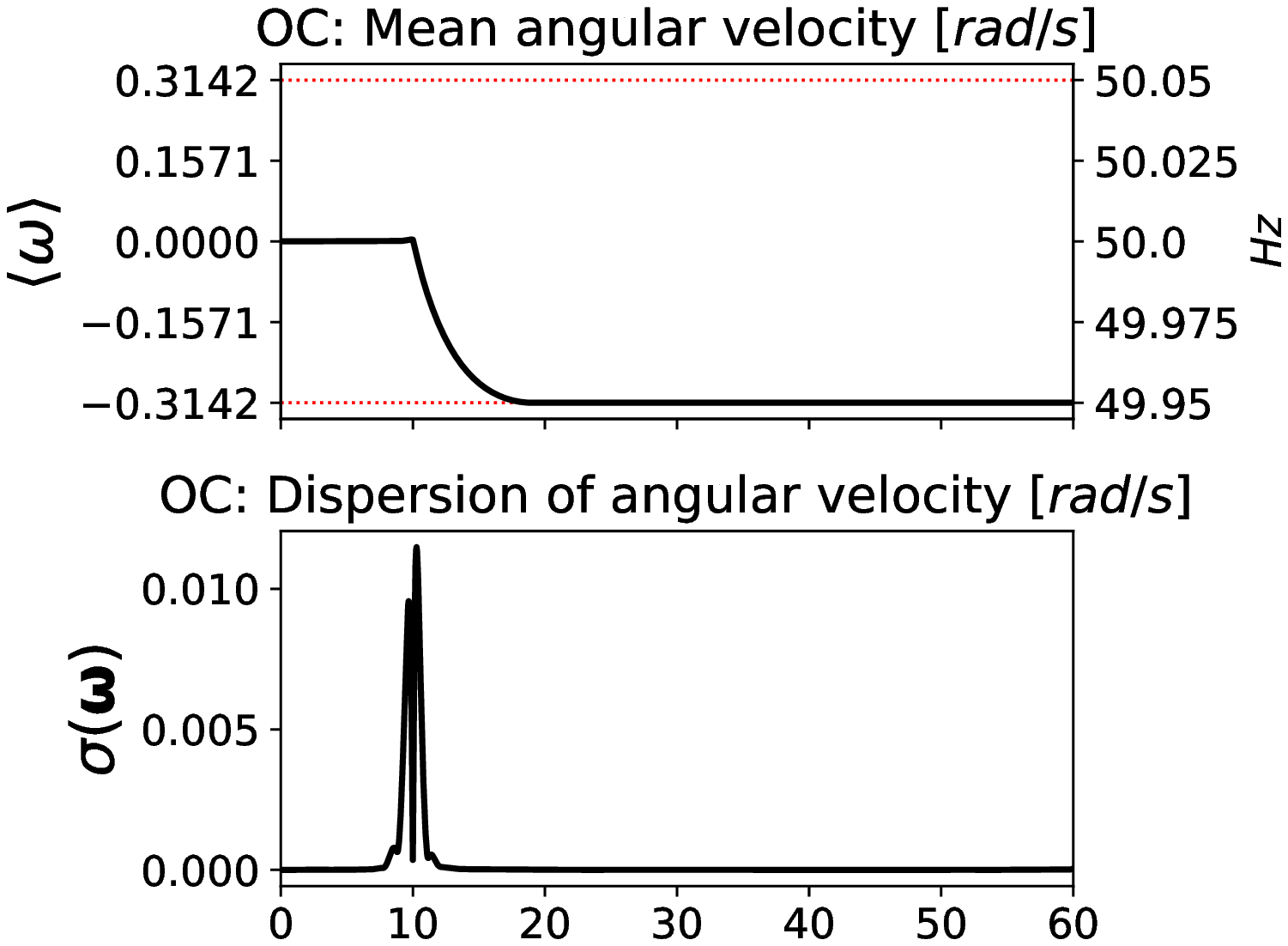}
		\label{fig:OC_persistent_system_stability}
	}
	\caption{Angular velocity mean and deviation in the test system under the persistent disturbance. Red dotted lines show operational limits. Angular velocity mean and deviation in the test system under the temporary disturbance. Red dotted lines show operational limits. Each control keeps the mean angular velocity $\langle \omega \rangle$ within its bounds and gradually synchronizes the system after each change in power by the disturbance. Notice that OC synchronizes the angular velocities to the boundary of its admissible set of values.}
	\label{fig:Persistent-System-Synchronization}
\end{figure*}

\end{document}